\documentclass[journal, 12pt, draftclsnofoot, onecolumn]{IEEEtran}
\usepackage[cmex10]{amsmath}
\usepackage{cite,array,datetime,calc}
\usepackage{caption,subcaption}
\usepackage{multicol,multirow}
\usepackage{ dsfont }

\renewcommand\tablename{Tab.}
\renewcommand\thetable{\arabic{table}}

\makeatletter\g@addto@macro\normalsize{
    \setlength\abovedisplayskip{3pt plus 3pt minus 2pt}
    \setlength\belowdisplayskip{3pt plus 3pt minus 2pt}
    \setlength\abovedisplayshortskip{0pt}
    \setlength\belowdisplayshortskip{10pt}
}\makeatother

\usepackage{pgfplots,tikz}
\pgfplotsset{compat=newest}
\newlength\figwidth
\newlength\figheight

\usepackage{flcond-s}

\usepackage{bm}
\usepackage{amssymb}
\usepackage{amsthm}

\global\long\def\Re{\mathbb{R}}

\global\long\def\0{\mathbf{0}}
\global\long\def\1{\mathbf{1}}
\global\long\def\A{\mathbf{A}}
\global\long\def\B{\mathbf{B}}
\global\long\def\a{\mathbf{a}}
\renewcommand{\b}{\mathbf{b}}
\renewcommand{\c}{\mathbf{c}}

\global\long\def\D{\mathcal{D}}
\global\long\def\E{\mathbf{E}}
\global\long\def\f{\mathbf{f}}
\global\long\def\fx{\f(\x)}
\global\long\def\fxa{\f(\xa)}
\global\long\def\F{\mathcal{F}}
\global\long\def\fo{\mathbf{f}_{0}}
\global\long\def\fox{\fo(\x)}
\global\long\def\g{\mathbf{g}}

\global\long\def\i{\mathcal{I}}
\global\long\def\I{\mathbf{I}}
\global\long\def\M{\mathbf{M}}
\global\long\def\p{\mathbf{p}}

\global\long\def\Rn{\Re^{n}}

\global\long\def\t{\mathbf{t}}
\global\long\def\u{\mathbf{u}}

\global\long\def\x{\mathbf{x}}
\global\long\def\cX{\mathcal{X}}
\global\long\def\ox{(\x)}
\global\long\def\oxh{(\hat\x)}
\global\long\def\xa{\mathbf{x^{\star}}}
\global\long\def\xh{\hat{\x}}
\global\long\def\y{\mathbf{y}}
\global\long\def\z{\mathbf{z}}

\global\long\def\la{\bm{\lambda}}

\newcommand{\muv}{\bm\mu}

\global\long\def\and{\text{and}}

\global\long\def\8{\ensuremath\infty}
\theoremstyle{plain}
\newtheorem{theorem}{Theorem}
\newtheorem{lemma}[theorem]{Lemma}
\newtheorem{proposition}[theorem]{Proposition}
\newtheorem{result}[theorem]{Result}
\newtheorem{corollary}[theorem]{Corollary}

\theoremstyle{definition}
\newtheorem{definition}[theorem]{Definition}
\newtheorem{remark}[theorem]{Remark}
\newtheorem*{remark*}{Remark}

\newtheorem*{condo}{{\color{blue}Old} Qualifying Conditions}

\newcommand{\q}{\mathfrak{q}}

\newcommand{\dd}{\,\mathrm{d}}
\newcommand{\de}{\mathrm{d}}

\newcommand{\cmin}{c_{\min}}
\newcommand{\cmax}{c_{\max}}

\newlength{\norlen}
\newcommand{\norm}[1]{\left|\hspace{-\norlen}\left|\hspace{-\norlen}\left|#1\right|\hspace{-\norlen}\right|\hspace{-\norlen}\right|}

\newcommand{\normi}[1]{\norm{#1}_\infty}
\newcommand{\vnorm}[1]{\left\|#1\right\|}

\newcommand{\iseti}{\mathbb A}
\newcommand{\isete}{\mathbb B}

\newcommand{\Rmp}{\ensuremath{\Re_+^m}}

\newcommand{\ssum}{\textstyle{\sum}}

\newenvironment{TempProof}{\begin{IEEEproof}}{\end{IEEEproof}}

\newcommand\New[1]{#1}

\let\qcold\qc
\renewcommand\qc[1]{\New{\qcold{#1}}}
\let\qccold\qcc
\renewcommand\qcc[1]{\New{\qccold{#1}}}

\begin{document}
\author{Martin~Jakobsson, 
            Carlo Fischione, 
            and~Pradeep~Chathuranga~Weeraddana
\thanks{The authors are with the Department of Automatic Control and the ACCESS Linnaeus Center, %
KTH Royal Institute of Technology, Stockholm, Sweden. Email: \{mjakobss, carlofi, chatw\}@kth.se.}%
\thanks{This work was supported by the EU projects NoE HYCON2 and STREP Hydrobionets.
}
\thanks{A preliminary version of this work was presented in~\cite{Jakobsson2012}.}}

\title{Extensions of Fast-Lipschitz Optimization}

\maketitle

\begin{abstract}
The need of fast distributed solvers for optimization problems in networked systems has motivated the recent development of the Fast-Lipschitz optimization framework. In such an optimization, problems satisfying certain qualifying conditions, such as monotonicity of the objective function and contractivity of the constraints, have a unique optimal solution obtained via fast distributed algorithms that compute the fixed point of the constraints. This paper extends the set of problems for which the Fast-Lipschitz framework applies. Existing assumptions on the problem form are relaxed and new and generalized qualifying conditions are established by novel results based on Lagrangian duality. It is shown for which cases of more constraints than decision variables, and less constraints than decision variables Fast-Lipschitz optimization applies. New results are obtained by imposing non strict monotonicity of the objective functions. \New{The extended Fast-Lipschitz framework is illustrated by a number of examples, including a non-linear optimal control problem.}
\end{abstract}

\section{Introduction}
\New{
Fast-Lipschitz optimization is a recently proposed class or problems where the special structure ensures that the optimal point is the fixed point of the constraints.\footnote{\New{We give a brief technical summary of Fast-Lipschitz optimization in Section~\ref{sec:Background-Work}.}} These problems can therefore be solved through simple and decentralized algorithms, making the framework an interesting alternative for distributed applications over networks. The special structure of Fast-Lipschitz problems makes them useful also in centralized applications, where it is guaranteed that the solutions can be obtained through simple system of equations that do not use the traditional Lagrangian approach. }

In a general networked optimization problem, the network nodes or agents must coordinate their actions to optimize a network-wide objective function.
When information such as nodes' objectives, constraints and decision variables are distributed among the nodes and physically scattered across the network, or if the amount of information is too large to collect centrally, it can be impractical or even impossible to centrally compute the solution. For example, collecting information in one place might be too expensive if the network has limited communication resources, or it may be too slow if the solution is needed at the local nodes in real time. In these situations, fast distributed solution algorithms must be used.

Distributed optimization has a long history, and much of the recent developments build upon the work of Tsitsiklis\cite{Tsitsikilis,Tsitsiklis1986}. Several approaches exist for solving these problems, such as primal and dual decomposition methods.
In these methods, the primal or the dual problem is decomposed into local subproblems solved at the nodes. The subproblems are coordinated through a centralized master problem, which is usually solved by gradient or subgradient methods \cite{Boyd2007}. These methods have found many applications in network utility maximization, e.g., \cite{Kelly1998,Low1999,Palomar2007}.
Due to the slow convergence of subgradient algorithms, recent works have explored higher order methods. For example, \cite{Wei2011a} replaces the subgradient step with a Newton-like step.
An other decomposition approach which is faster and more robust than the standard decomposition methods is the alternating direction method of multipliers (ADMM). The method has recently attracted a substantial interest, especially for problems with large data sets \cite{Gabay1976,Eckstein92,Schizas2008,Schizas2008b,Boyd2011}.

Although the decomposition methods mentioned above distribute the computational workload over the network nodes, the subproblems must still be coordinated among the nodes. For example, dual decomposition based methods must update dual variables in a central master problem. This requires the network to iteratively 1) transfer information form all nodes to some central ``master" node; 2) centrally update the dual variables; 3) broadcast the updated dual variables to all nodes of the network and back to 1) until convergence.

To avoid the centralized master node, peer-to-peer or multi-agent methods have recently been proposed to coordinate the subproblem through local neighbor interactions based on consensus algorithms, e.g., ~\cite{Jadbabaie2003,Boyd2006,Dimakis2008,Olshevsky2009,Chiuso2011}. In these algorithms, nodes update their decision variables as convex combinations of the decision variables of their neighbors without a centralized master.
In \cite{Nedic2009}, the consensus algorithm has been combined with gradient descent to solve an unconstrained optimization problem where the objective is a sum of local, convex functions. The work is extended in \cite{Nedic2010,SundharRam2010}, which investigate constraints and randomness. While the previous papers solve the primal problem, \cite{Duchi12,Tsianos12, Zhu2012} use consensus-based algorithms for the dual problem. Higher order methods are considered also for peer-to-peer optimization, e.g.,~\cite{zanella_et_al__2011__newton_raphson_consensus} solves an unconstrained primal problem, whereas~\cite{Zargham12} solves a linearly constrained dual problem, both by approximating Newton's algorithm through consensus.

Consensus based methods have many benefits, such as resilience to node failures and changing network topology. However, since every round of consensus requires  message passings, also these methods may suffer from communication overheads.
A recent study of the tradeoff between communication
and local computation can be found in \cite{Tsianos12b}.
The communication overhead is a problem especially in large scale distributed networks or wireless sensor networks, where the energy expenditure for communication can be orders of magnitude larger than the energy for computation~\cite{Fischione11}.

The methods discussed thus far assume convex problems. There are other classes of algorithms that do not rely necessarily on convexity, but on other structural properties. Three such classes are abstract optimization \cite{Notarstefano2009}, which generalizes linear optimization, monotonic optimization \cite{Tuy2000,Jorswieck2010,Bjornson2012}, where the monotonicity of the objective function is used to iteratively refine a solution within bounds of the feasible region, and Interference Function optimization~\cite{Zander,Grandhi94,Foschini1993,Yates95,Herdtner2000}, which is the fundamental framework to solve radio power control problems over wireless networks. Given the importance of Interference Function optimization as a precursor of Fast Lipschitz optimization, and considering that we will give some application examples later on in this paper, we give below some technical detail on such an optimization framework.

In \emph{Interference Function optimization}, typical radio power control problems are solved in a simple distributed (as in decentralized) way. In such an optimization approach, it is assumed that the nodes of a wireless network transmit signals of a certain level of radio power, say node~$i$ transmits with power $p_i$. The signal is received at the intended receiver corrupted by multiplicative wireless channel attenuations and additive interference by other transmitter nodes. The level of power that the transmit power of node $i$ has to overcome at the receiver in order to get the signal decoded is usually denoted by $\i_i(\p)$, the \emph{interference function} of  transmitter~$i$~\cite{Yates95}, where $\p=[p_{1}, p_{2},\ldots,p_{n}]^T$ is the vector of all transmit radio powers of the $n$ transmitters in the wireless network.
The goal of the basic power control problem is to  minimize the radio powers $p_{i}$, while overcoming the interference at each receiver, i.e.,
\begin{align}
\min_{\p} & \quad \p \nonumber \\
{\rm s.t.} & \quad p_{i}\ge\i_{i}(\p)\quad\forall i.\label{eq power problem}
\end{align}
\New{
Note that the problem above is a vector optimization, where the minimization is carried out with respect to the non-negative orthant (see Section~\ref{sub: vec opt} for details). Roughly speaking, the minimization of problem \eqref{eq power problem} makes all the components of $\p$ become small simultaneously.
The solution of problem~\eqref{eq power problem}} is a particularly successful instance of distributed optimization. Affine versions of $\i_i(\p)$ are the simplest and best studied type of interference function, but in theory one can consider functions $\i(\p)$ of any form. The first distributed algorithm to solve such a problem was proposed in \cite{Zander}, and improved in~\cite{Grandhi94,Foschini1993}. The algorithm was later generalized to the Interference Function framework by Yates~\cite{Yates95}.
In this framework, a function $\i(\p)$ is called \emph{standard} if, for all $\p\ge\0$, it fulfils
\begin{itemize}
\item \emph{Monotonicity}: If $\p\ge\p'$, then $\i(\p)\ge\i(\p')$,
\item \emph{Scalability}: For all $\alpha>1,\,\alpha\i(\p)>\i(\alpha\p)$.
\end{itemize}
When problem \eqref{eq power problem} above is feasible, and the functions $\i_{i}(\p)$ are standard, the unique optimal solution is given by the fixed point of the iteration
\begin{align} \label{iterations}
p_{i}^{k+1}:=\i_{i}(\p^{k}),
\end{align}
or $\p^{k+1}:=\i(\p^{k})$ in vector form. Here, $p_i^k$ is the power of transmitter $i$ at time $k$ and $\p^{k}=[p_{1}^{k}, p_{2}^{k},\ldots,p_{n}^{k}]^T$. The computation of the optimal solution by these iterations is much simpler than using the classical parallelization and decomposition methods. This is because there is no longer a need to centrally collect, compute and redistribute the coupling variables of the problem since $\i_{i}(\p^{k})$ can be know locally at node $i$~\cite{Yates95}. Even in a centralized setting, iteration \eqref{iterations} is simpler than traditional Lagrangian methods, since no dual variables need to be stored and manipulated.
The iterations require only that every node successively updates its transmit power by using local knowledge of other nodes' current decision variables (radio powers). Another advantage is that the algorithm converges even though such a knowledge is delayed, i.e., when the decision variables $p_j^k$ of other nodes come with some delay at node $i$ \cite{Bertsekas1997}.

Extensions of the Interference Function framework have been proposed in \cite{Sung2005}, where optimization problems whose interference functions are not standard by Yates' definition are investigated. Instead, they introduce Type-II standard functions, which for all $\p\ge\0$ fulfill
\begin{itemize}
  \item \emph{Type-II monotonicity}: If $\p\le\p',$ then $\i(\p)\ge\i(\p')$.
  \item \emph{Type-II scalability: }$\forall\alpha>1$, $\i(\alpha\p)>(1/\alpha)\i(\p)$.
\end{itemize}
These functions are shown to have the same fixed point properties as Yates' standard functions, i.e., problem \eqref{eq power problem} with Type-II standard constraints can be solved through repeated iterations of the constraints \eqref{iterations}.

Fast-Lipschitz optimization is a natural generalization of the Interference Function approach on how to solve distributed optimization problems over networks by using fixed point iterations similar to \eqref{iterations}, but when the constraints are neither standard nor Type-II standard~\cite{Fischione11}. It also considers more general objective functions $\g_0(\p)$ for problem \eqref{eq power problem}.
The framework considers a class of possibly non-convex and multi-objective problems with monotonic objective functions. These problems have unique Pareto optimal solutions, well defined by a contractive system of equations formed by the problem constraints. Therefore, Fast-Lipschitz problems are solved without having to introduce Lagrangian functions and dual variables or consensus based iterations. This makes the framework particularly well suited when highly decentralized solutions, with few coordination messages, are required. This is important in typical areas such as  wireless sensor networks and in multi-agent systems.

In this paper, we substantially extend the class of problems that are currently solvable with the Fast-Lipschitz framework.
\New{
In particular, we 1) introduce a new qualifying condition that unifies and extends the existing conditions of \cite{Fischione11} , 2) consider problems with more, or less, constraints than variables, which was not considered in~\cite{Fischione11}, and 3) study objective functions that are not strictly monotonic, which was not considered in~\cite{Fischione11}.

The remainder of this paper is organized as follows.
In Section \ref{sec:Background-Work} we clarify our notation and give a brief overview of the current state of Fast-Lipschitz optimization.
Section \ref{sec:Main} presents a new, general qualifying condition for Fast-Lipschitz optimization. The new qualifying condition is proved (using an approach different from \cite{Fischione11}) in~\ref{sub: main proof}.
In Section~\ref{sub: special cases} we give special cases of the general qualifying condition that have much less analytical and computational complexity, and highlight the connection of these special cases to existing qualifying conditions and related work.
Furthermore, Section~\ref{sec: relaxed FL form} relaxes some of the requirements of Fast-Lipschitz form by considering problems with more constraints than variables in Subsection~\ref{sub: additional}, and problems with less constraints than variables in~\ref{sub:Fewer-constraints-than}, while Subsection~\ref{sub: Non-strictly-monotonic-cost} refines some of the requirements on the objective function.
}
Section \ref{sub ex simple} features an example that illustrates some of the new results of this paper. \New{In Section~\ref{sub cex} the new results developed in this paper are applied to a family of optimal control problems, where the problem structure is utilized to determine the optimal solution without computations.
}
Finally, the paper is concluded in Section~\ref{sec:Conclusions-and-future}.

\section{Preliminaries}\label{sec:Background-Work}
This section clarifies notation and recalls Fast-Lipschitz optimization to provide the essential background definitions for the core contribution of this paper.

\subsection{Notation}
Vectors and matrices are denoted by bold lower and upper case letters, respectively. The components of a vector $\x$ are denoted $x_{i}$  or $ [\x]_{i}$. Similarly, the elements of the matrix $\A$ are denoted $A_{ij}$ or  $\left[\A\right]_{ij}$.
The transpose of a vector or matrix is denoted $\cdot^{T}$. $\I$  and $\1$  denote the identity matrix and the vector of ones. A vector or matrix where all elements are zero is denoted by $\0$.

The gradient $\nabla\fx$ \New{is the transpose of the Jacobian matrix, i.e.,}  $\left[\nabla\fx\right]_{ij}=\partial f_{j}\ox/\partial x_{i}$, whereas $\nabla_i\fx$ denotes the $i$th row of $\nabla\fx$. Note that $\nabla\fx^k=\left(\nabla\fx\right)^k$, which has not to be confused with the $k$th derivative.
The spectral radius of $\A$ is denoted $\rho(\A)$.
Vector norms are denoted $\vnorm{\cdot}$ and matrix norms are denoted $\norm{\cdot}$. Unless specified $\vnorm{\cdot}$ and $\norm{\cdot}$ denote arbitrary norms. $\normi{\A}=\max_{i}\sum_{j}|A_{ij}|$ is the norm induced by the $\ell_\infty$ vector norm.

All inequalities are intended \emph{element-wise}, i.e., they have nothing to do with positive definiteness.

\New{
\subsection{Vector optimization and Pareto optimality}\label{sub: vec opt}
In this paper we are concerned with maximization of vector valued objective functions $\fox\in\Re^m$. A vector optimization problem involves a proper cone $\mathcal K$ (see, e.g., \cite[Section 4.7]{Boyd2004}). In this paper we focus on the case when $\mathcal K$ is non-negative orthant
\[
\mathcal K = \Rmp\triangleq\{\y \, : \, y_i\ge0,\,\forall i=1,\dots,m\}.
\]
This maximization of a vector is formally expressed as
\begin{equation}\label{eq vecop 1}
\text{maximize (with respect to } \Rmp)\;  \fox.
\footnote{\New{
In this paper we only consider optimization with respect to the cone \Rmp. Therefore, for notational simplicity we simply write $\max\;\fox$, and the vector optimization~\eqref{eq vecop 1} should be understood whenever $\fox$ is vector valued.}
}
\end{equation}
The goal of the maximization~\eqref{eq vecop 1} is to find the decision vector x such that the components of $\fox$ are as big as possible with respect to the cone \Rmp.
In particular, when comparing two vectors $\x$ and $\y$ with respect to the cone \Rmp, we say $\x \succeq_{\Rmp} \y$ if $\x-\y \in \Rmp$ and $\x \preceq_{\Rmp}\y$ if $\y-\x\in\Rmp$.
 Note that this corresponds exactly to the component-wise inequalities ($\y\ge\x$ and $\y\le\x$), and for this reason we will use $\ge$ rather than $\succeq_{\Rmp}$ throughout the rest of the paper.
Unlike scalars, where it must hold that either $a\ge b$, or $a<b$, two vectors might not be comparable.
For example, for $\x=[1\;2]^T$ and $\y=[3\;1]^T$, we have $\x\not\ge\y$ and $\y\not\ge\x$.

Based on the discussion above, $\x_1$ is preferable to $\x_2$ in problem~\eqref{eq vecop 1} if $\fo(\x_1)\ge\fo(\x_2)$.
A feasible decision variable $\hat\x$ is said to be \emph{Pareto optimal} if there is no other feasible vector $\x$ such that $\fox\ge\fo(\hat\x)$. A problem can have several Pareto optimal points --  in this case each Pareto optimal point is incomparable to any other Pareto optimal points, but preferable to any point that is not Pareto optimal. If a problem only has one Pareto optimal point $\xa$, then $\xa$ is the unambiguously best choice among the decision vectors.

\emph{Scalarization} is a useful technique for generating Pareto optimal points. Scalarization is performed by picking a weight vector $\muv$ in the interior of $\Rmp$, i.e., $\muv >\0$, and solving the scalar optimization problem
\begin{equation}\label{eq vecop 2}
\max \; \muv^T \fox =\sum_{i=1}^m \mu_i [\fox]_i.
\end{equation}
Any point $\hat\x$ that is optimal in problem~\eqref{eq vecop 2} is Pareto optimal in~\eqref{eq vecop 1}.
Although scalarization does not necessarily generate all Pareto optimal points, one can show that if a point $\xa$ is optimal for \emph{all} scalarization vectors $\muv>\0$, then $\xa$ is the unique Pareto optimal point.
}

\subsection{Fast-Lipschitz optimization}
\label{sub: old main}
We will now give a formal definition of Fast-Lipschitz problems. For a thorough discussion of Fast-Lipschitz properties we refer the reader to the above mentioned paper.
\begin{definition}
A problem is said to be on\emph{ Fast-Lipschitz form} if it can be written
\begin{equation}
\label{eq FLform}
\begin{array}{cl}
\max & \fo(\x)\\
\text{s.t.} & x_{i}\le f_{i}(\x)\quad  \forall i\in\iseti \\
 & x_{i}=f_{i}(\x)\quad  \forall i\in\isete, \\
\end{array}
\end{equation}
where \begin{itemize}\setlength{\itemsep}{1ex}

  \item
  $\fo:\Rn\rightarrow\Re^{m}$ is a differentiable scalar ($m=1$) or vector valued ($m\ge2$) function,

  \item
  $\iseti$ and $\isete$ are complementary subsets of $\{1,\dots,n\}$,

  \item $f_{i}\,:\,\Rn\to\Re$ are differentiable functions.
  \end{itemize}
We will refer to problem \eqref{eq FLform} as our main problem. From the individual constraint functions we form the vector valued function $\f:\Rn\rightarrow\Rn$ as
  \(
  \f(\x)=\begin{bmatrix}f_{1}(\x) & \cdots & f_{n}(\x)\end{bmatrix}^{T}.
  \)
\end{definition}

\begin{remark}
The characteristic feature of Fast-Lipschitz form is a pairing such that each variable $x_i$ is associated to one constraint $f_i(\x)$. The form $x\le f(\x)$ is general, since any constraint on canonical form, $g(\x)\le0$, can be written $x\le x-\gamma g(\x)$ for some positive constant $\gamma$.
\end{remark}

\begin{definition}
For the rest of the paper, we will restrict our attention to a \emph{bounding box}
\(
\D=\left\{ \x\in\Rn\,|\,\a\le\x\le\b\right\}.
\)
We assume $\D$ contains all candidates for optimality and that $\f$ maps $\D$ into $\D$,
\(
\f:\D\to\D.
\)
This box arise naturally in practice, since any real-world quantity or decision must be bounded.
\end{definition}

\begin{definition}
A problem is said to be \emph{Fast-Lipschitz} when it can be written \New{in} Fast-Lipschitz form and admits a unique Pareto optimal solution $\xa$, defined as the unique solution to the system of equations
\begin{equation}
\xa=\f(\xa).\label{eq x=00003Df(x)}
\end{equation}
\end{definition}

A problem written  \New{in}  Fast-Lipschitz form is not necessarily Fast-Lipschitz. The following \emph{qualifying conditions} guarantee that a problem  \New{in}  Fast-Lipschitz form is Fast-Lipschitz.

\begin{condo}
\label{qual cond 0}
For all $\x$ in $\D$, $\fo(\x)$ and $\f(\x)$ should be everywhere differentiable and fulfill at least one of the following cases (e.g., (0) and (i), or (0) and (ii)):

\let\oldarraystretch\arraystretch  \let\oldtabcolsep\tabcolsep \newlength\qcw \newlength\ics
\newcommand{\mpe}[1]{
\begin{minipage}{\qcw}
\vspace*{-15pt}
{\begin{flalign} & #1 & \end{flalign}}
\vspace*{-22pt}
\end{minipage}}

\let\myline\relax

\setlength{\tabcolsep}{0pt}
\renewcommand{\arraystretch}{1.1} 
\setlength\qcw{130mm} 
\setlength\ics{7pt} 

\noindent\hspace*{\fill}
\begin{subequations}\label{all original QC}\begin{tabular}{r @{~} p{2cm} p{\qcw}}
 \myline
&\textbf{(0)} &  \mpe{\label{old 0} \nabla\fo(\x)>\0}\\[\ics] \myline

AND
&\textbf{(i.a)} &  \mpe{\label{old i.a} \nabla \fx\ge\0}\\\myline
&\textbf{(i.b)} &  \mpe{\label{old i.b} \norm{\nabla\f(\x)}<1}\\[\ics]\myline

OR
&\textbf{(ii.a)} &  \mpe{\color{blue}\label{old ii.a} \fo(\x)=c\1^T\x,\text{ with }c>0} \\\myline
&\textbf{(ii.b)} &  \mpe{\color{blue}\label{old ii.b} \nabla\f(\x)\le\0 \text{ (or more generally,  }\nabla\f(\x)^2\ge\0)} \\\myline
&\textbf{(ii.c)} &  \mpe{\label{old ii.c} \normi{\nabla\f(\x)}<1}\\[\ics]\myline

OR
&\textbf{(iii.a)} &  \mpe{\color{blue} \fo(\x)\in\Re\label{old iii.a}} \\\myline
&\textbf{(iii.b)} &  \mpe{\color{blue}\label{olddeltas} \normi{\nabla\f(\x)}<\dfrac{\bar\delta}{\bar\delta+\bar\Delta},
\quad\text{where}\quad
\begin{matrix}{
\bar\delta \triangleq { \min_{i}\min_{\x\in\D}}\nabla_{i}\fo(\x)},&\\\myline
\bar\Delta \triangleq { \max_{i}\max_{\x\in\D}}\nabla_{i}\fo(\x)
\end{matrix}}
\\[2pt] \myline
\end{tabular}\end{subequations}
\setlength{\tabcolsep}{\oldtabcolsep}
\renewcommand{\arraystretch}{\oldarraystretch}
\end{condo}

\begin{theorem}[\hspace{.01pt}{\cite[Theorem 3.3]{Fischione11}}]
\label{theorem: old main}
 A problem  \New{in}  Fast-Lipschitz form that fulfills any pair of the \New{Old} Qualifying Conditions is Fast-Lipschitz, i.e., it has a unique Pareto optimal point given by $\xa=\f(\xa).$
\end{theorem}

Once it is known that a problem is Fast-Lipschitz, computing the solution becomes a matter of solving the system of equations \eqref{eq x=00003Df(x)}, which in general is much easier than solving an optimization problem using Lagrangian multipliers. This is particularly evident when $\f(\x)$ is contractive on $\D$, a property assured by the qualifying conditions. In this case, the iterations $\x^{k+1}:=\f(\x^{k})$ converge geometrically to the optimal point $\xa$, starting from any initial point $\x^0\in\D$.

Fast-Lipschitz optimization problems need not be convex, but convex problems that fulfill the qualifying conditions can be rewritten and solved with the Fast-Lipschitz framework. For example, this is true for any problem where the constraints are  standard:

\begin{proposition}[{\hspace{-0.01ex}\cite[Theorems 4.2]{Fischione10}}]
  If problem \eqref{eq power problem} is feasible, and the constraints standard, then the problem is Fast-Lipschitz.
\end{proposition}

This concludes the preliminary part of the paper. We are now ready to present the core contributions of this paper.

\section{A general qualifying condition}
\label{sec:Main}
\New{This section presents a new qualifying condition, that generalizes and unifies the conditions of Section~\ref{qual cond 0}. The new condition is introduced in Subsection~\ref{sec: newCond}, together Theorem~\ref{theorem: Main theorem - New}, which formally states the role of the condition. Subsection~\ref{sub: algorithms} shows that problems fulfilling the general qualifying condition of Subsection~\ref{sec: newCond} have contractive constraints. These results enable finding the optimal point through fixed point iterations, and they also function as a preliminary result to the proof of Theorem~\ref{theorem: Main theorem - New} in Subsection~\ref{sub: main proof}.  }

From this point and onwards we will do small change of terminology. We will call each set of related conditions a ``Qualifying Condition", rather than a ``case" of the qualifying conditions. We will still use the notion of ``case" when referring to the Old Qualifying Conditions~\eqref{all original QC}.
For example, ``case (i)" will refer to the groups \eqref{old 0}-\eqref{old i.b} of the old qualifying conditions, while ``case (ii)" refers to  \eqref{old 0} and \eqref{old ii.a}-\eqref{old ii.c}.

\subsection{New Qualifying Conditions}
\label{sec: newCond}

\newlength\lblen
\newlength\qcondlen
\newlength\deflen

Consider once again the optimization problem $\eqref{eq FLform}$ on Fast-Lipschitz form, surrounded by the bounding box $\D$. Just as in Section~\ref{qual cond 0}, the qualifying conditions presented in this section are used to ensure that a problem  \New{in}  Fast-Lipschitz form also is Fast-Lipschitz.
\New{
In the upcoming qualifying conditions, we will use the ratio
\begin{equation}
\label{eq q}
 \q\ox \triangleq \min_j \frac{\min_i [\nabla\fox]_{ij}}{\max_i [\nabla\fox]_{ij}}\,.
\end{equation}
The value of $\q\ox$ is the smallest ratio of any two elements from the same column of $\nabla\fox$. When $\nabla\fox\ge\0$, $\q(\x)$ is always non-negative.
Furthermore, $\q(\x)\le 1$ by construction, with equality if and only if all rows of $\nabla\fox$ are identical. In fact, $\q\ox$ can be seen as a penalty for when the objective function gradient points in a different direction than the vector $\1$.

We are now ready to state the new qualifying conditions. We will begin with the most general form of the qualifying conditions. Special cases of this condition will be discussed in Section~\ref{sub: special cases}.

\New{
\small
\setcounter{qcn}{-1}
\renewcommand{\tablename}{Tab.}
\renewcommand{\thetable}{\arabic{table}}
\setlength\lblen{25mm}
\setlength\qcondlen{100mm}
\setlength\deflen{100mm}
\begin{center}
\begin{tabular}{|c|l|}
\hline
\multicolumn{2}{|c|}{\textbf{General Qualifying Condition}}
\\ \hline
 \textbf{\qc{Qk}} &
\begin{Tqcond}
\itemsep=1.2\itemsep
  \label{Qk}
  \case $\nabla\fo(\x)\ge\0$ with non-zero rows
  \label{Qk.f0}
  \case
  \label{Qk.cont}
  $\norm{\nabla\f(\x)}<1$ 
  \case[There exists a $k \in \{1,2,\dots\}\cup\infty$ such that ]
  \label{Qk.pos}
 \begin{tabular}[]{l}
  When $k<\infty$, then \\
  $ \nabla\f(\x)^k\ge\0$
 \end{tabular}
  \case
  \label{Qk.inf}
  \begin{tabular}[]{l}
  When $k>1$, then \\
  $  \normi{\sum_{l=1}^{k-1} \nabla\fx^l} <  \q  (\x)
  \triangleq \min_j \displaystyle\frac{\min_i [\nabla\fox]_{ij}}{\max_i [\nabla\fox]_{ij}}$
  \end{tabular}
\end{Tqcond}
 \\ \hline
\end{tabular}
\end{center}
}
In the condition above, we allow the parameter $k$ to be any positive integer, or infinity.
Each value of $k$ corresponds to a different case in the proof of Theorem~\ref{theorem: Main theorem - New} below. In the extremes $k=\infty$ and $k=1$, there is no need to fulfill conditions  \qcc{Qk.pos} and \qcc{Qk.inf} respectively. Both these cases will be discussed further in Section  \ref{sub: special cases}.}

We now give Theorem~\ref{theorem: Main theorem - New}, which plays a role analogous to that of Theorem~\ref{theorem: old main} in Section~\ref{sub: old main}.

\begin{theorem}
  \label{theorem: Main theorem - New}
  Assume problem \eqref{eq FLform} is feasible, and that \New{\qc{Qk}} holds for every $\x\in \D$. Then, the problem is Fast-Lipschitz, i.e., the unique Pareto optimal solution is given by
  $\xa=\f(\xa).$
\end{theorem}

\begin{TempProof}
The theorem is proved in Section~\ref{sub: main proof}.
\end{TempProof}

\begin{remark}
  It should be emphasized that \New{all qualifying conditions of this paper} are only sufficient, i.e., a problem that does not fulfil the qualifying conditions can still be Fast-Lipschitz.
  For example, by considering certain transformations of the constraint functions, it is possible to relax condition the condition \qcc{Qk.inf} that contains the norm $\normi{\cdot}$.
  One can also relax certain requirements on the Fast-Lipschitz form of problem $\eqref{eq FLform}$, (e.g. by considering problems with more variables than constraints or problem where the objective function only depends on a subset of the variables (see Section~\ref{sec: relaxed FL form}).
\end{remark}

\subsection{Contraction properties of Fast-Lipschitz problems}\label{sub: algorithms}
\New{In this subsection we briefly discuss the contractiveness of $\fx$. This is important since it allows a Fast-Lipschitz optimization problem to be solved through fixed point methods. These results will also be used in the proof of Theorem~\ref{theorem: Main theorem - New} in the next subsection.}

Once optimization problem \eqref{eq FLform} is shown to be Fast-Lipschitz, solving it becomes a matter of finding the  point $\xa=\fxa$. In a centralized setting, one can use any suitable method for solving such a system of equations, e.g., Newton like methods. However, if $\fx$ is contractive, the simplest way to solve $\xa=\fxa$ is to repeatedly evaluate $\fx$. This method works both in a centralized and a distributed setting. Compared to other distributed methods, it has the benefit of being totally decentralized, i.e., there is no master problem or coordinating \New{node. Furthermore, the iterations converge even when some of the nodes, due to dropped or delayed packets, only have access to outdated information of the neighbors' decisions, see \cite[Proposition 3.6]{Fischione11} or \cite{Bertsekas1997} for details.} We show below that if the \New{general qualifying condition~\qc{Qk} holds}, then $\fx$ is contractive.

\begin{proposition}[Sec. 3.1 in \cite{Bertsekas1997}]\label{lem: contractions}
  Let $\f$ be a mapping from a closed subset of $\Rn$ into itself, $\f:\cX\to\cX$. If there is a norm $\vnorm{.}$ and a scalar $\alpha<1$ such that
  \(
  \vnorm{\fx-\f(\y)} \le \alpha \vnorm{\x-\y}
  \)
  for all $\x,\y\in\cX$, then $\fx$ is a contraction mapping. As a result:
  \begin{itemize}
    \item $\xa=\fxa$ is the unique fixed point of $\f$ in $\cX$.
    \item For every initial point $\x^0\in\cX$, the sequence
    \(
    \x^{k+1}:=\f(\x^k)
    \)
     converges linearly to $\xa$.
  \end{itemize}
\end{proposition}
Since we know $\f:\D\to\D$, $\f$ is a contraction mapping if we can find a vector norm such that $\vnorm{\fx-\f(\y)} \le \alpha \vnorm{\x-\y}$ for all $\x,\y\in\D$.

\begin{lemma}\label{lem: f contrac}
  \New{If the general qualifying condition~\qc{Qk} holds,} $\fx$ is contractive.
\end{lemma}
\begin{TempProof}
  Parameterize the line between $\x,\y\in\D$ by
  \(
  \u(t)=t\x+(1-t)\y, \; 0\le t \le 1,
  \)
  and let $\g(t)=\f(\u(t))$. Then,
  \(
  \de \g(t)/\de t=\nabla\f\left(\u(t)\right)^T(\x-\y),
  \)
  wherefore
  \begin{align*}
  \fx-\f(\y)
  &=
  \g(1)-\g(0)=\int^1_0 \frac{\de\g(t)}{\de t} \dd t
  =\int^1_0\nabla\f\left(\u(t)\right)^T \de t \, (\x-\y)
  \triangleq
  \A (\x-\y),
  \end{align*}
  where each element in $\A$ is the integral of the corresponding element in $\nabla\f(\u(t))$.
  \New{Let $\norm{\cdot}_a$ be the matrix norm that satisfies condition~\qcc{Qk.cont}, and define $\norm{\cdot}_b$ such that $\norm{\A}_b \triangleq \norm{\A^T}_a$. It is straight-forward to show that  $\norm{\cdot}_b$ inherits the matrix norm properties of  $\norm{\cdot}_a$.\footnote{\New{See the appendix for a short proof.}}
  We can now bound
  \begin{align*}
  \norm{\A}_b &=
  \norm{ \int^1_0\nabla\f\left(\u(t)\right)^T \dd t }_b
  \le
  \int^1_0 \norm{\nabla\f\left(\u(t)\right)^T}_b \dd t
  \\&
  \le \int_0^1 \max_{\z\in\D} \norm{\nabla\f(\z)^T}_b \dd t = \max_{\z\in\D} \norm{\nabla\f(\z)^T}_b = \max_{\z\in\D} \norm{\nabla\f(\z)}_a
  \triangleq \alpha.
  \end{align*}
  The first inequality above is the triangle inequality. The second inequality holds since $\D$ is convex whereby $\u(t)\in\D$ for all $t$. The maximum exists since $\D$ is compact, and $\alpha<1$ by \qcc{Qk.cont}.}

  From \cite[Theorem 5.6.26]{Horn85}, we know that there exists an \emph{induced} matrix norm $\norm{.}_V$ such that $\norm{\M}_V \le \norm{\M}_a$ for every matrix $\M\in\Re^{n\times n}$. Let $\vnorm{.}_v$ be the vector norm that induces $\norm{.}_V$.
  By the properties of induced matrix norms we get
  \begin{align*}
  \vnorm{\fx-\f(\y)}_v &= \vnorm{\A (\x-\y)}_v \le \norm{\A}_V\vnorm{(\x-\y)}_v
  \le \norm{\A}_a\vnorm{(\x-\y)}_v \le \alpha\vnorm{(\x-\y)}_v
  \end{align*}
  as desired. Since $\fx:\D\to\D$, $\fx$ is a contraction mapping. This concludes the proof.
\end{TempProof}

We are now ready for the main proof of the paper.

\subsection{Proof of Theorem \ref{theorem: Main theorem - New}}\label{sub: main proof}
In this section we prove Theorem~\ref{theorem: Main theorem - New}, which is one of the main contributions of this paper. The proof will be given as a series of lemmas as outlined in the following.

\begin{TempProof}[Proof of Theorem \ref{theorem: Main theorem - New}]
Consider optimization problem \eqref{eq FLform}.
When the \New{general qualifying condition~\qc{Qk} holds} for all $\x\in\D$, then the following steps ensure that $\xa=\fxa$ is the unique optimal solution.
\begin{enumerate}
  \item First, we restrict ourselves to optimization problems on Fast-Lipschitz form with only inequality constraints, without loss of generality by Lemma \ref{lem: All inequalities} below.
  \item The inequality-only constrained optimization problem allows us to show that all feasible points of the optimization problem are regular~\cite{BertsekasNonlinear}, wherefore any optimal point $\xh$ must fulfill the KKT-conditions, see Lemma \ref{lem: regular} below.
  \item Any point that fulfills the KKT-conditions must be a fixed point of $\fx$, see Lemma \ref{lem: lambda > 0}.
  \item Finally, we show that there exists a unique fixed point $\xa=\fxa$ by Proposition~\ref{lem: contractions} and Lemma~\ref{lem: f contrac}. Therefore, $\xa$ is the only point fulfilling the KKT-conditions and $\xa$ must be the optimum.
\end{enumerate}
\end{TempProof}

The first lemma allows us to focus on problems consisting only of inequality constraints. This is an important step that permits us to use the KKT conditions (see, e.g., \cite[3.1.1]{Bertsekas1997}) to establish the existence an uniqueness of optimal solutions.
\begin{lemma}
\label{lem: All inequalities} If the inequality-only constrained optimization problem

\begin{align}
\label{eq inequality only}
\begin{array}{cl}
    \max & \fo(\x)\\
    \textnormal{s.t.} & x_{i}\le f_{i}(\x)\quad\forall i\in\{1,\dots,n\}  \qquad\qquad\qquad\qquad
\end{array}
\\
\intertext{is Fast-Lipschitz, then so is any problem}
\label{eq mixed relations}
\begin{array}{cl}
    \max & \fo(\x)\\
    \textnormal{s.t.} & x_{i}\le f_{i}(\x)\quad\forall i\in\iseti       \phantom{,}\quad\qquad  \iseti \cup \isete=\{1,\dots,n\}   \\
     & x_{i}=f_{i}(\x)\quad\forall i\in\isete,     \quad\qquad \iseti \cap \isete=\O
\end{array}
\end{align}
obtained by switching any number of the inequalities for equalities.
\end{lemma}

\begin{TempProof}
  Let $\F_{\ref{eq inequality only}}$ and
  $\F_{\ref{eq mixed relations}}$ be the feasible regions of
  problem~\eqref{eq inequality only} and
  problem~\eqref{eq mixed relations} respectively. The point $\xa=\f(\xa)$ is feasible in both problems
  and by definition uniquely optimal for
  problem~\eqref{eq inequality only} since it is Fast-Lipschitz.
  Suppose, contrary to the lemma, that
  problem~\eqref{eq mixed relations} is not Fast-Lipschitz.
  Then there exists some feasible point
  $\xh\in\F_{\ref{eq mixed relations}}\subset
  \F_{\ref{eq inequality only}}$ such
  that $\fo(\xh)\ge\fo(\xa)$ which contradicts the unique optimality
  of $\xa$ in problem~\eqref{eq inequality only}.
\end{TempProof}

With Lemma \ref{lem: All inequalities} in mind, we develop the rest of the proof of Theorem \ref{theorem: Main theorem - New} by focusing on the inequality-only constrained problem~\eqref{eq inequality only} instead of main problem~\eqref{eq FLform}. The inequality-only problem on canonical form is
\begin{equation}\label{eq x-f<0}
  \begin{array}{cl}
\min & -\fo(\x)\\
\text{s.t.} & g_i(\x) = x_i - f_i\ox \le 0\quad\forall i.
\end{array}
\end{equation}

\begin{definition}
In problem \eqref{eq x-f<0}, a point $\x$ is \emph{regular} if the gradients of all active constraints at $\x$ form a linearly independent set (see \cite{BertsekasNonlinear}).
\end{definition}

\begin{lemma}
  \label{lem: regular}
  If problem \eqref{eq inequality only} fulfills \New{the general qualifying condition~\qc{Qk}}, then every point $\x\in\D$ in problem \eqref{eq x-f<0} is regular.
\end{lemma}
\begin{TempProof}
The gradients of the individual constraints $g_i(\x)$ are the columns of $\nabla\g(\x)=\I-\nabla\f(\x)$.
Since \New{condition \qcc{Qk.cont} implies
\(
\rho\big(\nabla\f(\x)\big)\le\norm{\nabla\f(\x)}<1,
\)}
the eigenvalues of $\nabla\g(\x)$ lie in a ball of radius $\rho\big(\nabla\f(\x)\big)<1$, centered at $1$. Hence, no eigenvalue of $\nabla\g(\x)$ is zero and $\nabla\g(\x)$ is invertible, wherefore the constraint gradients $\nabla\g_i(\x)$ (the columns of $\nabla\g\ox$) form a  linearly independent set.
\end{TempProof}

\New{We now scalarize  problem~\eqref{eq x-f<0}, by considering}
\begin{equation}\label{eq scale}
\begin{array}{cl}
\min & -\muv^{T}\fo(\x)\\
\text{s.t.} & g_i(\x) = x_{i} - f_{i}(\x) \le 0 \quad\forall i.
\end{array}
\end{equation}
\New{for a positive vector $\muv\in\Re^m$ (see Subsection~\ref{sub: vec opt}).}
Since an arbitrary scaling of $\bm\mu$ does not change the solution of the problem, we restrict ourselves to $\bm\mu$ fulfilling $\sum_k\mu_k=1$.

Introduce dual variables $\la\in\Rn$,  and form the Lagrangian function
\(
L(\x,\la)=-\muv^{T}\fo(\x)+\la^{T}(\x-\f(\x)).
\)
Since every $\x\in\D$ in problem~\eqref{eq scale} is regular (Lemma~\ref{lem: regular}), any pair $(\hat{\x},\hat{\la})$ of locally optimal variables must satisfy the KKT-conditions (see e.g., \cite[3.1.1]{Bertsekas1997}). In particular,
$\hat{\x}$ must be a minimizer of $L(\x,\hat{\la})$, which requires
\begin{equation}\label{eq nabla_L}
 \nabla_{\x}L(\hat{\x},\hat{\la})=
-\nabla\fo(\hat{\x})\muv+\hat{\la}-\nabla\f(\hat{\x})\hat{\la}=\0,
\end{equation}
and complementarity must hold, i.e.,
  \begin{equation}
  \label{eq strict_comp}
  \hat{\lambda}_{i}(\hat{x}_{i}-f_{i}(\hat{\x}))=0\quad\forall i.
  \end{equation}

We will soon show that \qc{Qk} implies \mbox{$\hat\la>\0$}.
\New{To this end, the following remark is useful.
\begin{remark}
  \label{lem: non-neg}
   Let $\la=\A\c$. If $\c>\0$ and $\A\ge\0$ with non-zero rows, then $\la>\0$.
\end{remark}
\noindent The statement above is trivial, but we give it as Remark~\ref{lem: non-neg}  since we will refer to it several times throughout the paper. Note that $\A\ge\0$ and $\c>\0$ is not sufficient for $\A\c>\0$, this is the reason that condition~\qcc{Qk.f0} requires non-zero rows.
}

The following lemma is the main part of the proof of Theorem \ref{theorem: Main theorem - New}, and establishes that the optimal dual variable is strictly positive.
\begin{lemma}\label{lem: lambda > 0}
Whenever  \New{\qc{Qk} holds}, any pair $(\hat{\x},\hat{\la})$ with $\hat\x\in\D$ satisfying  equation~\eqref{eq nabla_L} must have $\hat\la>\0$.
\end{lemma}
\begin{TempProof}
For notational convenience, let us fix one
$\x\in\D$, which we denote $\xh$, and introduce
\begin{equation}
\label{eq A def}
\A \triangleq \nabla\f(\xh)\in\Re^{n\times n}
\quad
\and\quad\c \triangleq \nabla\fo(\xh)\muv\in\Rn.
\end{equation}
Note that condition \New{\qcc{Qk.f0}} and Remark~\ref{lem: non-neg} give $\c>\0$ for every $\bm{\mu}>\0$.

With the new notation, equation \eqref{eq nabla_L} can be written as
\begin{equation}\label{eq basic}
-\c+\hat{\la}-\A\hat{\la}=\0,\,\text{ or }\,\hat{\la}=(\I-\A)^{-1}\c
\end{equation}
whenever $\I-\A=\I-\nabla\f(\hat\x)$ is invertible, which is true for all $\x\in\D$ (proof of Lemma~\ref{lem: regular}). The expansion of this inverse gives
\New{\begin{subequations}
\label{eq lambda main proof}
\begin{align}
\label{eq infinite series}
\hat{\la}&=(\I+\A+\A^{2}+\dots)\c
\\
\nonumber &= (\I+\A^{k}+\A^{2k}+\dots)(\I+\A+\dots+\A^{k-1})\c
\\
\label{eq q proof 0} &=
\underbrace{\ssum_{l=0}^\infty (\A^{k})^l}_\text{``left matrix"} \;\underbrace{ (\I+\A+\dots+\A^{k-1}) \c}_\text{``right vector"}.
\end{align}
\end{subequations}}
\New{The first step is showing that the ``left matrix" above is non-negative with non-zero rows. If $k\to\infty$, the ``right vector" becomes identical to the right hand side of~\eqref{eq infinite series}. The ``left matrix" must therefore equals identity\footnote{\New{A different way to see this is by noting that $\lim_{k\to\infty}\A^k=0$ since $\rho(|A)\1$ by condition~\qcc{Qk.cont}. The ``left matrix" therefore evaluates to $\sum_{l=1}^\infty \0^l$, where only the term $\0^0=\I$ gives a contribution.}}, which is non-negative with non-zero rows. For all other $k<\infty$, the ``left matrix" is ensured non-negative by condition \qcc{Qk.pos} (with non-zero rows guaranteed by the first term $\A^0=\I$).

When the ``left matrix" is non-negative with non-zero rows, a sufficient condition for $\hat\la>\0$ is that the right vector is positive (Remark~\ref{lem: non-neg}), i.e.,}
\begin{equation}
\label{eq q proof 1}
(\I+\underbrace{\A+\dots+\A^{k-1}}_{\triangleq \B} )\c > \0 \quad \Leftrightarrow \quad -\B \c < \c.
\end{equation}
\New{When $k=1$, then $\B=\0$ and \eqref{eq q proof 1} holds trivially, since $\c>\0$. For $k>1$,} let $\cmin$ and $\cmax$ be the minimum and maximum elements of $\c$ and consider row $i$ of inequality \eqref{eq q proof 1}.
We can now bound the right side by $\cmin \le c_i$ and the left side by
\begin{align}
\label{eq Bc bound}
\Big|[-\B\c]_i\Big| & =\left|\ssum_{j=1}^n B_{ij}c_j\right|
\le \ssum_{j=1}^n |B_{ij}||c_j|
\le \max_i\ssum_{j=1}^n |B_{ij}|\cmax
= \normi{\B}\cmax.
\end{align}
Therefore, equation \eqref{eq q proof 1} holds  if $\normi{\B}\cmax<\cmin$, or
\begin{equation}\label{eq q proof 2}
\normi{\B} < \frac\cmin\cmax.
\end{equation}
Let $\nabla_i\fo\oxh$ be the $i$th row of $\nabla\fo(\xh)$ and define
\[
\a(\muv) = \underset{\nabla_i\fo\oxh}{\text{argmin }}{\nabla_i\fo\oxh  \muv },
 \quad\text{ and }\quad
\b(\muv) = \underset{\nabla_i\fo\oxh}{\text{argmax }}{ \nabla_i\fo\oxh  \muv }.
\]
Let $d_k(\muv)=a_k(\muv)/b_k(\muv)$ and express the components of $\a$ as $a_k(\muv)=d_k(\muv)b_k(\muv)$.
Since \New{$\c=\nabla\fo\oxh\muv$}, we have
 \[
 \cmin=\a(\muv)\nabla\fo\oxh = \sum_k a_k(\muv) \mu_k = \sum_kd_k(\muv)b_k(\muv) \mu_k
 \]
and
  \[
  \cmax = \b(\muv)\nabla\fo\oxh = \sum_k b_k(\muv)\mu_k.
  \]
The fraction in \eqref{eq q proof 2} can therefore be bounded by\vspace{2pt}
\begin{align}
\frac{\cmin}{\cmax}
&
=\frac{\sum_k a_k(\muv) \cdot \mu_k}{\sum_k b_k(\muv) \cdot \mu_k}
= \frac{\sum_k d_k(\muv) \cdot b_k(\muv) \cdot \mu_k}{\sum_k b_k(\muv) \cdot \mu_k}
\label{eq q proof 3}
\ge \frac{\sum_k d_{\min}(\muv) \cdot b_k(\muv) \cdot \mu_k}{\sum_k b_k(\muv) \cdot \mu_k}=d_{\min}(\muv),
\end{align}
\vspace*{3pt}
where
\begin{align}
d_{\min}(\muv)
&
=\min_k d_k(\muv)
=\min_k\frac{a_k(\muv)}{b_k(\muv)}
\label{eq q proof 4}
\ge \min_k\frac{\min_i[\nabla\fo(\xh)]_{ik}}{\max_i[\nabla\fo(\xh)]_{ik}}
=\q(\xh)
\end{align}
with $\q\ox$ defined in \eqref{eq q}. By the definition of $\B$ and condition \qcc{Qk.inf} we get
\[
\normi{\B}
=\normi{\sum_{l=1}^{k-1} \nabla\f(\xh)^{\New{l}}}
< \q(\xh),
\]
which together with inequalities \eqref{eq q proof 2}-\eqref{eq q proof 4} ensures
\[
\normi{\B}
=\normi{\ssum_{l=1}^{k-1} \nabla\f(\xh)}
< \q(\xh)
\le d_{\min}
\le \frac{\cmin}{\cmax},
\]
wherefore $\hat\la>\0$ by inequalities \eqref{eq q proof 1}-\eqref{eq q proof 2} and Remark~\ref{lem: non-neg}.
\end{TempProof}
We now know that every pair $(\xh,\hat\la)$ satisfying the KKT conditions must have $\hat\la>\0$, provided that \qc{Qk} holds. Furthermore, strict complementarity (equation \eqref{eq strict_comp}) must hold. Since $\hat\lambda_i>0$, we have
\mbox{$\hat\lambda_i(\hat x_i-f_i\oxh)=0$}
only if \mbox{$\hat x_i-f_i\oxh=0$}, i.e., $\xh=\f(\xh)$. In other words, any candidate for primal optimality must be a fixed point of $\fx$.
It remains to show that there always exists a fixed point, and that this fixed point is unique. But this is already done, since \qc{Qk} and Lemma~\ref{lem: f contrac} ensure $\fx$ is contractive, and the result follow from Proposition~\ref{lem: contractions}.

By now, we have shown that the unique point $\xa=\fxa$ is the only possible optimum of problem \eqref{eq scale}. Since this is true for \New{\emph{any} scalarization vector $\muv>\0$}, $\xa$ is the unique Pareto optimal point of problem \eqref{eq inequality only}. Finally, Lemma~\ref{lem: All inequalities} extends this to the originally considered problem~\eqref{eq FLform}. By this, we have taken all the steps to prove Theorem~\ref{theorem: Main theorem - New}.\hfill$\square$

In the next section we revisit the \qc{Qk}, through a number of special cases.

\section{Special cases of \qc{Qk}} \label{sub: special cases}
\New{
In this section, we discuss the general qualifying condition (\qc{Qk}) in more detail. Moreover, we present several new qualifying conditions, each of which implies~\qc{Qk}.

\qc{Qk} has the benefit of giving a unified view of the qualifying conditions. This is convenient for proving properties of Fast-Lipschitz problems, and also for giving an overall understanding for what the qualifying conditions ensure.
However, \qc{Qk} may not always be suitable  when determining whether or not a given problem (or class of problems) is Fast-Lipschitz. This is because the generality of \qc{Qk} comes at the price of analytical and computational complexity  and cumbersome notation.
For example, conditions \qcc{Qk.pos} and \qcc{Qk.inf} become increasingly tedious to verify as the integer $k$ grows.
However, as we will see, the special cases can yield clean and easily verifiable conditions, which are much easier to use in practice than \qc{Qk}.

Furthermore, the specialized cases provide easy comparison to the old qualifying conditions, and other related work such as the standard and type-II standard function of \cite{Yates95} and \cite{Sung2005}.
We start with the simplest case \qc{Q1}, given below.
}

\setlength\lblen{30mm}
\setlength\qcondlen{90mm}
\setlength\deflen{90mm}
\begin{center}
\small
\New{
\begin{tabular}{|c|l|}
\hline
\multicolumn{2}{|c|}{\textbf{Qualifying Condition \ref{Q1}}}
\\ \hline
 \textbf{\qc{Q1}} &
\begin{Tqcond}
  \label{Q1}
  \case
  \label{Q1.f0}
  $\nabla\fo(\x)\ge\0$ with non-zero rows
  \case
  \label{Q1.cont}
  $ \norm{\nabla\fx}<1$
 \case
  \label{Q1.pos}
  $\nabla\f(\x)\ge\0$
\end{Tqcond}
\\ \hline
\end{tabular}
}\end{center}

\New{
Qualifying condition \qc{Q1} is the special case of  \qc{Qk} when $k=1$.
It is the simplest case of Fast-Lipschitz optimization, and only requires a monotonic objective function $\fox$ and a monotonic, contractive constraint function $\fx$. \qc{Q1} is highly related to standard interference functions \cite{Yates95}. In fact, any problem~\eqref{eq FLform} with monotonic objective function and standard constraints is Fast-Lipschitz \cite[Theorem 4.2]{Fischione10},\cite{Jakobsson2014}.
The difference between \qc{Q1} and case~(i) of the Old Qualifying Conditions~\eqref{all original QC} lies in condition \qcc{Q1.f0}, where we now allow $\nabla\f_0(\x)\ge\0$ as long as no row consists only of zeros.
}

\New{

Qualifying Condition 2 is a simplified version of~\qc{Qk} with $k=2$.
\begin{center}
\small
\begin{tabular}{|c|l|}
\hline
\multicolumn{2}{|c|}{\textbf{Qualifying Condition \ref{Q2}}}
\\ \hline
 \textbf{\qc{Q2}} &
 \begin{Tqcond}
 \label{Q2}
  \case
  \label{Q2.f0}
  $\nabla\fo(\x)>\0$
  \case
  \label{Q2.pos}
  $\nabla\f(\x)^2\ge\0, \quad \big($e.g., $\nabla\f(\x)\le\0\big)$
  \case \label{Q2.inf}
   $ \displaystyle \normi{\nabla\f(\x)} < \q\ox \triangleq \min_j \frac{\min_i [\nabla\fox]_{ij}}{\max_i [\nabla\fox]_{ij}}$
\end{Tqcond}
\\ \hline
\end{tabular}
\end{center}
}

\New{
}
\New{
\begin{proposition}
 Qualifying condition \qc{Q2} implies \qc{Qk}
\end{proposition}
}
\New{
\begin{TempProof}
Condition \qcc{Qk.f0} is implied by \qcc{Q2.f0}.
Note that if any element of $\nabla\fox$ is zero, then $\q\ox=0$ and condition \qcc{Q2.inf} cannot be fulfilled. We can therefore, without loss of generality, use a strict inequality in \qcc{Q2.f0}.
Condition \qcc{Qk.cont} is fulfilled since \mbox{$\normi{\nabla\fx}<\q\ox$} by condition \qcc{Q2.inf}, and $\q\ox\le 1$.
Finally, conditions \qcc{Qk.pos} and \qcc{Qk.inf} are (for $k=2$) given exactly by \qcc{Q2.pos} and \qcc{Q2.inf}, respectively.
\end{TempProof}
}
\New{
Condition \qcc{Q2.pos} requires that the square of the gradient is non-negative, but is particularly easy to verify when $\nabla\fx\le\0$.
Note that also
$\nabla\fx\ge\0$ fulfills \qcc{Q2.pos}, i.e., \qcc{Q2.pos} is more general than \qcc{Q1.pos}.  However, this generalization comes at a cost since \qcc{Q2.inf} is more restrictive than \qcc{Q1.cont} in the sense that it requires the specific norm $\normi\cdot$ and that $\q\ox$ in general is less than one.

The formulation where $\nabla\fx\le\0$ corresponds to a non-increasing objective function, and the norm  $\normi{\nabla\fx}$ is small enough, \qc{Q2} is closely related to type-II standard functions \cite{Jakobsson2014}.
}

\New{
\renewcommand{\tablename}{Tab.}
\renewcommand{\thetable}{\arabic{table}}
\setlength\lblen{30mm}
\setlength\qcondlen{90mm}
\setlength\deflen{90mm}
\begin{center}\small
\begin{tabular}{|c|l|}
\hline
\multicolumn{2}{|c|}{\textbf{Qualifying Condition \ref{Qinf}}}
\\ \hline
 \textbf{\qc{Qinf}} &
 \begin{Tqcond}
 \label{Qinf}
  \case
  \label{Qinf.f0}
  $\nabla\fo(\x)>\0$

  \case
  \label{Qinf.inf}
  $ \displaystyle \normi{\nabla\f(\x)}< \frac{\q\ox}{1+\q\ox}$
\end{Tqcond}
\\ \hline
\end{tabular}
\end{center}
}
\New{
Qualifying condition \qc{Qinf} can be seen as the special case of \qc{Qk} when $k=\8$, as we see in Proposition~\ref{prop: Q3}.
In contrast to the other qualifying conditions, \qc{Qinf} does not require non-negativity of $\nabla\fx$ (or $(\nabla\fx)^k$).
It is only required that $\fo$ is monotonic and the infinity norm of $\nabla\fx$ is small enough.
However, \qcc{Qinf.inf} is significantly stricter than \qcc{Q2.inf} because $0 \le \q\ox \le 1$ implies $\normi{\nabla\fx} < \q\ox/(1+\q\ox)\le1/2$ in  \qcc{Qinf.inf}.
}
\New{
\begin{proposition}\label{prop: Q3}
  Qualifying condition \qc{Qinf} implies \qc{Qk}.
\end{proposition}
}
\New{
\begin{TempProof}
Condition \qcc{Qinf.f0} implies \qcc{Qk.f0} and  condition \qcc{Qinf.inf} implies \mbox{$\normi{\nabla\fx}<1/2$}, whereby condition \qcc{Qk.cont} is fulfilled.
Condition \qcc{Qk.pos} is irrelevant when $k=\infty$.%
\footnote{\New{
This is discussed in the proof of Theorem~\ref{theorem: Main theorem - New}, after equation~\eqref{eq lambda main proof}.
  }}
Finally, condition \qcc{Qk.inf} is implied by \qcc{Qinf.inf} because
\begin{equation}
\label{eq prop17 1}
\normi{\nabla\fx}<\frac{\q\ox}{1+\q\ox}
\; \implies \;
\frac{\normi{\nabla\fx}}{1- \normi{\nabla\fx}} < \q\ox,
\end{equation}
where the implication above follows from that
\[
 \normi{\ssum_{l=1}^{\infty} \nabla\f(\x)^l}
 \le \ssum_{l=1}^{\infty}\normi{ \nabla\f(\x)^l}
 \le \ssum_{l=1}^{\infty}\normi{\nabla\f(\x)}^l
 = \displaystyle\frac{\normi{\nabla\fx}}{1- \normi{\nabla\fx}}
 <\q\ox.
\]
The first two inequalities follow from sub-additive and sub-multiplicative properties of matrix norms. The equality follows from the geometric series because $\normi{\nabla\fx}<1$, and the last inequality is expression~\eqref{eq prop17 1}. This concludes the proof.
\end{TempProof}
}

\New{
The qualifying conditions can be further simplified by introducing
$\delta(\x) \triangleq \min_{i,j}\left[\nabla\fo(\x)\right]_{ij}$
and
$\Delta(\x) \triangleq \max_{i,j}\left[\nabla\fo(\x)\right]_{ij}$.
These are the smallest and largest elements of $\nabla\fox$, regardless of column. The difference compared to  $\bar\delta$ and $\bar\Delta$ in \eqref{olddeltas} is that we now evaluate $\delta$ and $\Delta$ for each $\x$, instead of taking the extremes over all $\x$. We can now bound
\begin{equation}\label{eq q-delta}
\q\ox
= \min_j \frac{\min_i [\nabla\fo]_{ij}}{\max_i [\nabla\fo]_{ij}}
\ge  
 \frac{\min_{ij} [\nabla\fo]_{ij}}{\max_{ij} [\nabla\fo]_{ij}}
= \frac{\delta\ox}{\Delta\ox}.
\end{equation}
}\New{
Since both $\q\ox$  and $\q\ox/(1+\q\ox)$ are increasing in $\q$ (recall $\q\ox\ge\0)$, we can lower bound $\q\ox$ by $\delta\ox/\Delta\ox$ in any one of the qualifying conditions above. This gives the remaining qualifying conditions of this section. They are all special cases of previous conditions --- easier to verify and analyze, at the expense of being more conservative.
}

\begin{table}[h]
\small
\New{
\renewcommand{\tablename}{Tab.}
\renewcommand{\thetable}{\arabic{table}}
\setlength\lblen{30mm}
\setlength\qcondlen{90mm}
\setlength\deflen{90mm}
\begin{center}
\begin{tabular}{|c|l|}
\hline
\multicolumn{2}{|c|}{\textbf{Qualifying Condition \ref{Q2d}}}
\\ \hline
 \textbf{\qc{Q2d}} &
 \begin{Tqcond}
 \label{Q2d}
  \case \label{Q2d.f0}
  $\nabla\fo(\x)>\0$

  \case \label{Q2d.pos}
  $\nabla\f(\x)^2\ge\0, \quad \big($e.g., $\nabla\f(\x)\le\0\big)$

  \case \label{Q2d.inf}
   $ \displaystyle \normi{\nabla\f(\x)} < \frac{\delta(\x)}{\Delta(\x)} $
\end{Tqcond}
\\ \hline
\multicolumn{2}{|c|}{\textbf{Qualifying Condition \ref{Qinfd}}}
\\ \hline
 \textbf{\qc{Qinfd}} &
 \begin{Tqcond}
 \label{Qinfd}
  \case \label{Qinfd.f0}
  $\nabla\fo(\x)>\0$

   \case \label{Qinfd.inf}
  $ \displaystyle \normi{\nabla\f(\x)} < \frac{\delta(\x)}{\delta(\x)+\Delta(\x)}$
  \vspace{-2pt}
\end{Tqcond}
\\ \hline
\end{tabular}
\end{center}
}
\end{table}
\New{
Qualifying conditions \qc{Q2d} and \qc{Qinfd} are obtained when inserting inequality~\eqref{eq q-delta} in \qcc{Q2.inf} and \qcc{Qinf.inf} respectively. They imply \qc{Q2} and \qc{Qinf} which in turn imply \qc{Qk} by construction. Note that qualifying conditions \qc{Q2d} and \qc{Qinfd} have previously appeared as~\cite[case~(ii)-(iii)]{Jakobsson2012}.
}

\New{
We will now end this section with the observation that the old qualifying conditions are a special case of \qc{Qk}.

\begin{proposition}
  The Old Qualifying Conditions~\eqref{all original QC} imply \qc{Qk}.
\end{proposition}
\begin{TempProof}
We will show how each case of the old qualifying conditions imply one of the qualifying conditions 1-5 above, which in turn implies \qc{Qk}.

\textbf{case~(i) $\Rightarrow$~\qc{Q1}}: Condition \qcc{Q1.f0} is implied by \eqref{old 0}, \qcc{Q1.cont} and \qcc{Q1.pos} are the same as \eqref{old i.a} and \eqref{old i.b}.

  \textbf{case~(ii) $\Rightarrow$ ~\qc{Q2d}}: Conditions \qcc{Q2d.f0} and \qcc{Q2d.pos} are the same as \eqref{old 0} and \eqref{old ii.b}. Condition~\eqref{old ii.a} implies $\nabla\fox=c\1$, whereby $\delta\ox=\Delta\ox=c$. Condition~\qcc{Q2d.inf}  therefore requires that \mbox{$\normi{\nabla\fx}<1$}, which is ensured by \eqref{old ii.c}.

  \textbf{case~(iii) $\Rightarrow$~\qc{Qinfd}}: Condition~\eqref{old iii.a} requires $\fox$ to be scalar valued, whereby $\nabla\fox$ only has one column. The deltas of condition~\eqref{olddeltas}  can therefore be written as
  $\bar\delta  \triangleq  \min_{\x\in\D}\delta\ox \le \delta\ox$
 and
 $\bar\Delta \triangleq \max_{\x\in\D}\Delta\ox \ge \Delta\ox$. This gives
 \begin{equation}
 \label{eq delta ineq}
\frac{\delta_{\min} }{\Delta_{\max}}
\le
\frac{\delta\ox }{\Delta\ox},
\end{equation}
whereby
\newcommand\deltafrac{\bar\delta/\bar\Delta}
\newcommand\deltafracx{\delta\ox /\Delta\ox}
\[
\normi{\nabla\fx}
<
\frac{\bar\delta }{\bar\delta +\bar\Delta }
=
\frac{\deltafrac}{1+\deltafrac}
\le
\frac{\deltafracx}{1+\deltafracx}
=
\frac{\delta\ox}{\delta\ox +\Delta\ox}
\]
and \qcc{Qinfd.inf} is ensured. The first inequality is condition~\eqref{olddeltas}, the second inequality follows from \eqref{eq delta ineq} because $h(a)=a/(1+a)$ is an increasing function of $a$.

We have now showed how each case of the Old Qualifying Conditions implies \qc{Qk}, through the implication chain
\newcommand{\implied}{\rotatebox[origin=c]{180}{$\implies$}}
\[
\begin{tabular}{c c c c c c l}
        & \raisebox{-5pt}{\rotatebox[origin=c]{25}{\implied}} & \qc{Q1} & \implied & (i) & \\

    \qc{Qk} & \implied & \qc{Q2} & \implied & \qc{Q2d} & \implied & (ii) \\

        & \raisebox{7pt}{\rotatebox[origin=c]{-25}{\implied}} & \qc{Qinf} & \implied & \qc{Qinfd}   & \implied & (iii).
\end{tabular}
\]
This concludes the proof.
\end{TempProof}
The next section will loosen some of the assumptions of the optimization problem structure, i.e., investigate problems not entirely  \New{in}  Fast-Lipschitz form.
}

\New{
\section{Relaxations of the Fast-Lipschitz form}
\label{sec: relaxed FL form}
This section considers relaxations of the Fast-Lipschitz form that, e.g., require the same number of constraints as variables. Subsection~\ref{sub: additional} shows a technique for handling more constraints than variables, and Subsection~\ref{sub:Fewer-constraints-than} shows a situation with fewer constraints than variables. Finally, Subsection~\ref{sub: Non-strictly-monotonic-cost} discuss the case when the objective function does not depend on all variables.
}
\subsection{Additional constraints} \label{sub: additional}
In this section we  complement the Fast-Lipschitz form \eqref{eq FLform} with an additional set $\cX$. Hence, we consider the following problem:
\begin{equation}
  \label{eq FLform extended}
  \begin{array}{cl}
      \max & \fo(\x)\\
      \text{s.t.} & x_{i}\le f_{i}(\x)\quad  \forall i\in\iseti \\
       & x_{i}=f_{i}(\x)\quad  \forall i\in\isete \\
       & \x\in\cX.
  \end{array}
\end{equation}

\New{
\begin{corollary}
\label{prop: additional constraint}
  If \New{\qc{Qk}} holds, and $\cX$ contains a point
  $\xa=\fxa,$
  then problem \eqref{eq FLform extended} is Fast-Lipschitz.
\end{corollary}}
\begin{TempProof}
  Relax the problem by removing the new constraint $\x\in\cX$. The relaxed problem is our main problem \eqref{eq FLform}, whereby the qualifying conditions and Theorem~\ref{theorem: Main theorem - New} ensure $\xa=\fxa$ is the unique optimum. Since $\xa\in\cX$, this is also the unique optimum of problem~\eqref{eq FLform extended}.
\end{TempProof}

In theory, we can handle any set $\cX$ provided we can show $\xa\in\cX$. For example, $\cX$ does not need to be bounded, convex, or even connected (i.e., $\cX$ can consist of mutually disconnected subsets). In practice however, the most common form of $\cX$ is a box constraint, for example a requirement of non-negativity, $\cX=\{\x\,:\,\0\le\x\}$. In these cases, $\cX$ becomes the natural choice for the imagined bounding box $\D$.

\subsection{Fewer constraints than variables -- Constant constraints\label{sub:Fewer-constraints-than}}

The Fast-Lipschitz form in problem \eqref{eq FLform} requires one (and only one) constraint $f_i$ for each variable $x_i$.
In this section we will look at the case when the number of constraints ($f_i$) are fewer than the number of variables. We will assume that the individual variables are upper and lower bounded, which is always the case for problems of engineering interest, such as wireless networks.
This means we get an extra constraint set $\cX$ as discussed in Section \ref{sub: additional}.
We investigate the case when all constraints are inequalities. It follows from Lemma~\ref{lem: All inequalities} that the following results are true also for problems with equality constraints.

Consider a partitioned variable  $\x=\left[\y^{T}\,\z^{T}\right]^{T}\in\cX\subset\Rn$
and the problem
\begin{equation}\label{eq fewer constraints original-1}
\begin{array}{cl}
\max & \fo(\x)\\
\text{s.t.} & \y\le\f_{\y}(\x)\\
 & \x
\in\cX,
\end{array}
\quad\text{with}\quad
\cX=\left\{ \x\,:\,
\Bigg\{\hspace{-0.4em}
\begin{array}{l}
\y\in\cX_\y=  \{\y\,:\,\a_\y\le\y\le\b_\y\}
\\
\z\in\cX_\z= \{\z\,:\,\a_\z\le\z\le\b_\z\}
\end{array}
\right\}.
\end{equation}
In the formulation above, there are no constraints $f_i$ for the variables $z_i$.
However, by enforcing $\z\in\cX_\z$ twice we get the equivalent
problem
\begin{equation}
\begin{array}{cl}
\max & \fo(\x)\\
\text{s.t.} & \y\le\f_{\y}(\x)\\
 & \z\le\f_{\z}(\x)=\b_{\z}\\
 &\x\in\cX.
\end{array}\label{eq fewer constraints extended-1}
\end{equation}
 This problem has the right form \eqref{eq FLform extended} and is Fast-Lipschitz if
\[
\nabla\f
=\begin{bmatrix}
    \nabla_{\y}\f_{\x}\ox & \nabla_{\y}\f_{\z}\ox\\
    \nabla_{\z}\f_{\y}\ox & \nabla_{\z}\f_{\z}\ox
\end{bmatrix}\quad\text{and}\quad\nabla\fo=\begin{bmatrix}
\nabla_{\y}\fo\ox\\
\nabla_{\z}\fo\ox
\end{bmatrix}
\]
fulfills \New{\qc{Qk}}. Since $\f_{\z}\ox=\b_{\z}$
is constant, $\nabla\fx$ simplifies to
\[
\nabla\f=\begin{bmatrix}
\nabla_{\y}\f_{\y}\ox & \0\\
\nabla_{\z}\f_{\y}\ox & \0
\end{bmatrix}.
\]
the qualifying conditions. Moreover, they add nothing to either of  $||\cdot||_{1}$ or $||\cdot||_{\infty}$ and does therefore not have an impact on the qualifying conditions.
The special structure of $\nabla\fox$ can be exploited to construct less restrictive qualifying conditions.

Fist, consider problem \eqref{eq fewer constraints original-1} with a fixed $\z\in\cX_\z$. The problem can be written
\begin{equation}
  \label{eq subproblem}
  \begin{array}{cl}
    \max & \f_{0|\z}(\y)\\
    \text{s.t.} & \y\le\f_{\y|\z}(\y)\\
     & \y\in\cX_\y.
  \end{array}
\end{equation} We will refer to problem \eqref{eq subproblem} as the \emph{subproblem} $(\f_{0|\z}, \, \f_{\y|\z})$.

\begin{proposition}
\label{prop: fewer constraints}
Consider problem \eqref{eq fewer constraints extended-1}. Suppose that
(a) the subproblem  $(\f_{0|\z}, \, \f_{\y|\z})$
fulfills \qc{Qk} for all $\z\in\cX_\z$, and it holds for all $\x\in\cX$,  that~either
  \begin{enumerate}
  \setlength\topsep{0pt}
  \setlength\itemsep{0pt}
  \setlength\labelwidth{0pt}

  \item[\New{(b.i)}] \mbox{$\nabla_\z\fo(\x)\ge\0$,} and $\nabla_{\z}\f_{\y}(\x)\ge\0$ with non-zero rows, ~~or

  \item[\New{(b.ii)}] $\nabla_{\z}\fo(\x)\ge\0$ with non-zero rows, and $\nabla_\z\f_\y(\x)\ge\0$, ~~or

  \item[\New{(b.iii)}]\vspace{5pt}
  \(\displaystyle\frac{\normi{\nabla_{\z}\f_{\y}(\x)}}{1-\normi{\nabla_{\y}\f_{\y}(\x)}} <  \frac{\delta_{\z}(\x)}{\Delta_{\y}(\x)},\)
  where
  \(
  \delta_{\z}(\x)=\min_{ij}\left[\nabla_{\z}\fo(\x)\right]_{ij}
  \) and \(
  \Delta_{\y}(\x)=\max_{ij}\left[\nabla_{\y}\fo(\x)\right]_{ij}.
  \)\vspace{5pt}

  \end{enumerate}
Then,  problem \eqref{eq fewer constraints original-1} is Fast-Lipschitz.
\end{proposition}

\begin{TempProof}
We refer to the arguments of Section~\ref{sub: main proof} where we modify Lemma~\ref{lem: lambda > 0} as follows.

The particular form and partitioning  remains in the definitions \eqref{eq A def}
, giving
\begin{equation}
  \label{eq A}
  \A=\left[
  \begin{array}{cc}
    \nabla_{\y}\f_{\y}(\xh) & \0\\
    \nabla_{\z}\f_{\y}(\hat{\x}) & \0
  \end{array}
  \right]=\left[
  \begin{array}{cc}
    \A_{11} & \0\\
    \A_{21} & \0
\end{array}
\right]
 \quad\and\quad
\c=\begin{bmatrix}
\nabla_{\y}\fo(\xh)\\
\nabla_{\z}\fo(\xh)
\end{bmatrix}\muv
=\begin{bmatrix}
\c_{1}\\
\c_{2}
\end{bmatrix}.
\end{equation}
Consider again equation \eqref{eq basic} and denote $\E \triangleq (\I-\A)^{-1}$, i.e., $\la=\E\c$. As in the proof of Lemma~\ref{lem: lambda > 0}, $\E$ is well defined if $\rho(\A)<1.$ This is the case, since the eigenvalues of a block triangular matrix are the union of the eigenvalues of the diagonal blocks, wherefore
 $\rho(\A)=\rho(\A_{11})<1$
by assumption~(a).

As in the proof of Lemma \ref{lem: lambda > 0}, we must show $\la>\0$. This time we will make use of the block structure of $\A$ and $\c$.%
\footnote{Formulas for the inverse of a block matrix, as well as the products of two block matrices can be found in \cite{Matrix Cookbook}.}
From the block matrix inverse formula we get
\[
\la=(\I-\A)^{-1}\c
=\begin{bmatrix}
    \E_{11} & \E_{12}\\
    \E_{21} & \E_{22}
\end{bmatrix}\begin{bmatrix}
    \c_{1}\\
    \c_{2}
\end{bmatrix},
\]
 where $\E_{11}=\left(\I-\A_{11}\right)^{-1}$, $\E_{12}=\0$,
$\E_{21}=\A_{21}\E_{11}$ and $\E_{22}=\I$, i.e.,
\[
\begin{bmatrix}
    \la_1 \\
    \la_2
\end{bmatrix}
=\begin{bmatrix}
    \E_{11}\c_{1}\\
    \E_{21}\c_{1}+\c_{2}
\end{bmatrix}=\begin{bmatrix}
    \E_{11}\c_{1}\\
    \A_{21}\E_{11}\c_{1}+\c_{2}
\end{bmatrix}.
\]
Note that
$\la_{1}
=\E_{11}\c_{1}
=\left(\I-\A_{11}\right)^{-1}\c_{1},$
so $\la_1>\0$ since the subproblem $(\f_{0|\z},\f_{y|\z})$ fulfills Lemma~\ref{lem: lambda > 0}  by assumption (a).
The second block component is
\begin{equation}
\label{eq lambda_2}
\la_{2}=\A_{21}\E_{11}\c_{1}+\c_{2}=\A_{21}\la_{1}+\c_{2}.
\end{equation}
Given that $\la_{1}>\0$, we need to show $\la_2>\0$ if either of assumptions \New{(b.i)-(b.iii)} hold.

We start with assumption \New{(b.i)}, which ensures $\c_2\ge\0$ and $\A_{21}\ge\0$ with non-zero rows (so $\A_{21}\la_1>\0$ by Remark~\ref{lem: non-neg}), wherefore $\la_2>\0$ by equation \eqref{eq lambda_2}.

Assumption \New{(b.ii)} assures $\A_{21}\ge\0$ and $\c_2>\0$, wherefore $\la_2>\0$  by equation \eqref{eq lambda_2}.

Finally, assuming \New{(b.iii)} is fulfilled, we see that $\la_2>\0$ if
\(
\c_{2}>-\A_{21}\E_{11}\c_{1}.
\)
Analogous to equations \eqref{eq Bc bound} and \eqref{eq q proof 2}, this holds if
\(
\min_i[\c_{2}]_i>\normi{\A_{21}\E_{11}}\max_i[\c_{1}]_i,
\)
 or since $\c_1>\0$ and $\E_{11}=\left(\I-\A_{11}\right)^{-1}=\sum_{k=0}^{\infty}\A_{11}^{k}$,
when
\begin{equation}\label{eq t}
\frac{\min_i[\c_{2}]_i}{\max_i[\c_{1}]_i}>\normi{\A_{21}\ssum_{k=0}^{\infty}\A_{11}^{k}}.
\end{equation}
By the triangle inequality, the sub-multiplicative property of matrix norms and the geometric series, the right side of equation \eqref{eq t} can be upper bounded by
\begin{equation}
\label{eq t-2}
\frac{\normi{\A_{21}}}{1-\normi{\A_{11}}}
\ge
\normi{\A_{21}\ssum_{k=0}^{\infty}\A_{11}^{k}}
\end{equation}
Finally the definitions of $\c_i$ and $\sum_{j}\mu_{j}=1$ gives
\begin{align}
\frac{\min_i[\c_{2}]_i}{\max_i[\c_{1}]_i}
&
=\frac{\min_{i}\sum_{j}\left[\nabla_{\z}\fo(\x)\right]_{ij}\mu_{j}}
{\max_{i}\sum_{j}\left[\nabla_{\y}\fo(\x)\right]_{ij}\mu_{j}}
\label{eq t-1}
\ge \frac{\min_{ij}\left[\nabla_{\z}\fo(\x)\right]_{ij}\sum_{j}\mu_{j}}
{\max_{ij}\left[\nabla_{\y}\fo(\x)\right]_{ij}\sum_{j}\mu_{j}}
=  \frac{\delta_{\z}(\x)}{\Delta_{\y}(\x)}.
\end{align}
 By combining inequalities \eqref{eq t}-\eqref{eq t-1}, a sufficient condition
ensuring  $\la_{2}>\0$ is
\[
\frac{\normi{\A_{21}}}{1-\normi{\A_{11}}}<\frac{\delta_{\z}(\x)}{\Delta_{\y}(\x)},
\]
which is guaranteed by assumption (d). This concludes the proof.
\end{TempProof}

\subsection{Non-strictly monotonic objective function -- Variables missing in
objective function \label{sub: Non-strictly-monotonic-cost}}

Sometimes it is practical or necessary to formulate problems where
not all variables appear in the objective function.

For example, the problem
\begin{equation*}
  \begin{array}{c l}
    \max         &  \fo(\x)  \\
    \text{s.t.} &  \x  \le \f_x(\x,\z) \\
                &  \z \le \f_z(\x,\z)
  \end{array}
\end{equation*}
has some variables not affecting the objective function and is not  \New{in}  Fast-Lipschitz form. Redefining $\fo=\fo(\x,\z)$
gives a problem of the right form \eqref{eq FLform}, but $\nabla_{\z}\fo(\x,\z)=\0$
everywhere. Therefore, condition \qcc{Qk.f0} and Theorem \ref{theorem: Main theorem - New} can not be used to classify the problem as Fast-Lipschitz.

The situation above is a special case of the following problem.
Consider a partitioned optimization variable $(\x,\z)$ and the
problem
\begin{equation}
\begin{array}{cl}
\max & \fo(\x,\z)\\
\text{s.t.} & \x\le\f_{\x}(\x,\z)\\
 & \z\le\f_{\z}(\x,\z).
\end{array}\label{eq missing vars}
\end{equation}
Suppose $\fo$ is monotonic, i.e. $\nabla\fo\ge\0$, but partitioned such that
\[
\nabla\fo(\x,\z)=\begin{bmatrix}\nabla_{\x}\fo(\x,\z)\\
\nabla_{\z}\fo(\x,\z)
\end{bmatrix},
\]
where $\nabla_{\x}\fo(\x,\z)$ has non-zero rows for all $(\x,\z)\in\D$, while $\nabla_{\z}\fo(\x,\z)$ can have zero rows.
By partitioning the objective function gradient, one can find
situations when problem \eqref{eq missing vars} is actually Fast-Lipschitz.
\begin{proposition}
\label{prop: missing vars in obj}
Consider problem \eqref{eq missing vars}.
If it holds, for all $\x$ in $\D$, that
\begin{enumerate}
  \def\theenumi{(\alph{enumi}}
  \setlength{\itemsep}{1ex}
  \vspace{1ex}

  \item $\nabla_x \fo(\x,\z) > \0$ and $\nabla_z \fo(\x,\z) \ge \0$
  \item $\nabla\f(\x,\z)
  =\begin{bmatrix}
      \nabla_{\x}\f_{\x}(\x,\z) & \nabla_{\x}\f_{\z}(\x,\z)\\
      \nabla_{\z}\f_{\x}(\x,\z) & \nabla_{\z}\f_{\z}(\x,\z)
  \end{bmatrix}
  \ge\0,$

  \item \New{$\norm{\nabla\fx}<1$ for some matrix norm, and}

  \item $\nabla_{\z}\f_{\x}(\x,\z)$ has non-zero rows.
\end{enumerate}
Then, problem \eqref{eq missing vars} is Fast-Lipschitz.
\end{proposition}
\begin{remark}
The condition that the $i^{th}$ row of $\nabla_{\z}\f_{\x}(\x,\z)$ is
non-zero means that an increase in the variable $z_{i}$ will allow an increase of
some variable $x_{j}$, which in turn will influence the objective.
\end{remark}
\begin{TempProof}
The proof of Theorem~\ref{theorem: Main theorem - New} can be reused, with some alterations to Lemma~\ref{lem: lambda > 0}.

The partitioning of $\nabla\f$ and $\nabla\fo$  remains in $\A$ and $\c$, i.e.,
\[
\A
= \begin{bmatrix}
    \nabla_{\x}\f_{\x}(\x,\z) & \nabla_{\x}\f_{\z}(\x,\z)\\
    \nabla_{\z}\f_{\x}(\x,\z) & \nabla_{\z}\f_{\z}(\x,\z)
\end{bmatrix}
=\begin{bmatrix}
    \A_{11} & \A_{12}\\
    \A_{21} & \A_{22}
\end{bmatrix}
\quad\and\quad
\c
=\begin{bmatrix}\nabla_{\x}\fo(\x,\z)\\
\nabla_{\z}\fo(\x,\z)
\end{bmatrix}\muv
=\begin{bmatrix}
    \c_{1}\\
    \c_{2}
\end{bmatrix}.
\]

Assumption (a) and $\muv>\0$ gives $\c_{1}>\0$ and $\c_{2}\ge\0$.
Just as in Lemma~\ref{lem: lambda > 0}, assumptions (b) and (c) guarantee the existence and non-negativity of $\E=\left(\I-\A\right)^{-1}$, so  $\la=\left(\I-\A\right)^{-1}\c\ge\0$ is well defined and non-negative.%
\footnote{In contrast to Lemma~\ref{lem: lambda > 0}, we only have the weak inequality $\c\ge\0$.}
Thus, it remains to show $\la=\E\c>\0$.

Expressing the inverse $\E=\left(\I-\A\right)^{-1}$ block-wise, we have
\[
\left(\I-\A\right)^{-1}=\begin{bmatrix}\I-\A_{11} & \A_{12}\\
\A_{21} & \I-\A_{22}
\end{bmatrix}^{-1}=\begin{bmatrix}\E_{11} & \E_{12}\\
\E_{21} & \E_{22}
\end{bmatrix},
\]
 where
\(
\E_{11}=\left(\I-\left(\A_{11}+\A_{12}\left(\I-\A_{22}\right)^{-1}\A_{21}\right)\right)^{-1}
\)
 and
\(
\E_{21}=\left(\I-\A_{22}\right)^{-1}\A_{21}\E_{11}.
\)
We now have
\[
\begin{bmatrix}\la_{1}\\
\la_{2}
\end{bmatrix}=\begin{bmatrix}\E_{11} & \E_{12}\\
\E_{21} & \E_{22}
\end{bmatrix}\begin{bmatrix}\c_{1}\\
\c_{2}
\end{bmatrix}=\begin{bmatrix}\E_{11}\c_{1}\\
\E_{21}\c_{1}
\end{bmatrix}+\begin{bmatrix}\E_{12}\\
\E_{22}
\end{bmatrix}\c_{2}.
\]
Since $\E\ge\0$, all blocks $\E_{ij}$ are non-negative. The second term is always non-negative and can be ignored, it is enough to show that the first term is strictly positive.

Since $\c_{1}>\0$, and $\E_{11}\ge\0$, we have that $\la_{1}>\0$ if $\E_{11}$ has non-zero
rows (Remark~\ref{lem: non-neg}). This is always the case since $\E_{11}$ (defined as an inverse) is invertible.
The second component can be expressed in terms of the first component:
\(
\la_{2}  
= \left(\I-\A_{22}\right)^{-1}\A_{21}\la_{1}.
\)
When $\la_{1}>\0$, $\A_{21}\la_{1}>\0$ if $\A_{21}=\nabla_{\z}\f_{\x}(\x,\z)$
has non-zero rows (Remark~\ref{lem: non-neg}). This is true by assumption~(d).
Since $\left(\I-\A_{22}\right)^{-1}$ is invertible it has non-zero rows
and $\la_{2}=\left(\I-\A_{22}\right)^{-1}\left(\A_{21}\la_{1}\right)>\0$, which concludes the proof.
\end{TempProof}

We end this section by a general example.

\subsubsection*{Example}
Start \New{with} problem \eqref{eq FLform}. For a lighter notation, we assume all constraints are inequalities. Transforming the problem to an equivalent problem on 
 epigraph form gives
\[
\begin{array}{c l}
  \max & \t \\
  \text{s.t.} &
  \begin{array}[t]{l @{\,}l}
    \t&\le\fo(\x),\\
    \x&\le\f(\x).
  \end{array}
\end{array}
\]
This problem has a (non-strictly) monotonic objective, regardless of  $\fo(\x)$. By writing this as
\[
\begin{array}{c l @{}l}
  \max & \multicolumn{2}{l}{\g_0(\t,\x)=\t} \\
  \text{s.t.} &  \t &\le\g_t(\t,\x)=\fo(\x) \\
              &  \x &\le \g_x(\t,\x)=\f(\x)

\end{array}
\]
we obtain
\[
\nabla\g_0(\t,\x) =
\begin{bmatrix}\nabla_{\t}\g_0(\t,\x) \\
\nabla_{\x}\g_0(\t,\x)
\end{bmatrix} =
\begin{bmatrix}\I \\
\0
\end{bmatrix},\quad\and
\]
\[
\nabla\g(\t,\x)=\begin{bmatrix}\nabla_{\t}\g_{\t}(\t,\x) & \nabla_{\t}\g_{\x}(\t,\x)\\
\nabla_{\x}\g_{\t}(\t,\x) & \nabla_{\x}\g_{\x}(\t,\x)
\end{bmatrix}=
\begin{bmatrix} \0 & \0\\
\nabla\fo(\x) & \nabla\f(\x)
\end{bmatrix}.
\]

Proposition \ref{prop: missing vars in obj} can now be applied, and the problem is Fast-Lipschitz if $\nabla\g\ge\0$, $\rho(\nabla\g)<1$ and $\nabla\fo$ has non-zero rows. Since the eigenvalues of a block triangular matrix are the union of the eigenvalues of the diagonal blocks, we have $\rho(\nabla\g)=\rho(\nabla\f)$, wherefore the problem is Fast-Lipschitz if
\begin{itemize}
  \item $\nabla\fx\ge \0$, $\norm{\nabla\fx}<1$, and
  \item $\nabla\fox\ge\0$ with non-zero rows.
\end{itemize}
This is precisely qualifying condition~\qc{Q1} applied to the original problem \eqref{eq FLform}, i.e., no generalization was achieved by considering the epigraph form of the problem.

This concludes the main part of the paper. In the following section we will illustrate the new theory with two examples.

\section{Examples}\label{sec:Examples}
\New{We begin by illustrating Fast-Lipschitz optimization on a non-convex optimization example in Section~\ref{sub ex simple}. In Section \ref{sub cex}, we apply the novel results established in this paper to state conditions for when a non-linear optimal control problem is easily solvable.
}
\subsection{Simple non-convex example}
\label{sub ex simple}
Consider the problem
\begin{equation}
\label{eq toy example}
\begin{array}{cl}
\max & \fo(\x)\\
\text{s.t.} & \x \le \f(\x)\\
 & \x\in\cX=\{\x\,:\,\0\le\x\le\1\},
\end{array}
\end{equation}
where $\x\in\Re^2$,
\[
\fo(\x)=\begin{bmatrix} 2x_1 + x_2\\ x_1+2x_2 \end{bmatrix},
\quad\text{ and }\quad
\f(\x)=0.5\begin{bmatrix} 1+a x_2^2 \\1+bx_1^2 \end{bmatrix}.
\]
If either $a$ or $b$ are positive, the problem is not convex (the canonical constraint functions $x_i-f_i(\x)$ become concave).

At this point we do not know if $\xa=\fxa\in\cX$. However, following the results of Section~\ref{sub: additional}, we can assume this is the case and examine whether
\begin{equation}
\label{eq toy no add}
\begin{array}{cl}
\max & \fo(\x)\\
\text{s.t.} & \x \le \f(\x)
\end{array}
\end{equation}
fulfills \New{\qc{Qk} (or \qc{Q1}-\qc{Qinfd})}. The qualifying conditions must apply in the box $\D$, which we select equal to $\cX$ (no point outside of $\D=\cX$ is feasible, hence all feasible points lie in $\D$).

The  gradients of the objective and constraint functions are
\[
\nabla\fo=\begin{bmatrix}2&1\\1&2\end{bmatrix}
\quad\text{ and }\quad
\nabla\f=\begin{bmatrix}0&bx_1\\ax_2&0\end{bmatrix}.
\]
Since $\nabla\fo>\0$ for all $\x$, \New{ the assumptions  \qcc{Qk.f0}  on the objective function} is always fulfilled. We will now check all sign combinations of $a$ and $b$.

\paragraph*{$a,b\ge0$}

If both $a$ and $b$ are non-negative, then $\nabla\f(\x)\ge\0$ for all $\x$ in $\D$ and condition \qcc{Q1.pos} holds. To verify condition  \qcc{Q1.cont}, one must find a norm $\norm{\cdot}$ such that  $\norm{\nabla \f(\x)}<1$ for all $\x$ in $\D$. When the $\infty$-norm is used, we get
\[
\max_{\x\in\D} \normi{\nabla\f(\x)}=\max_{\x\in\D} \max\{|bx_1|,|ax_2|\} \le \max\{b,a\}<1
\]
if $a,b<1$. Thus, when $0\le a,b<1$, the problem is Fast-Lipschitz by \qc{Q1}.

\paragraph*{$a,b\le0$}

If $a$ and $b$ are instead non-positive, we have $\nabla\f(\x)\le\0$ for all $\x$ in $\D$ and condition \qcc{Q2d.pos} is fulfilled. In order to verify  \qcc{Q2d.inf}, we need $\delta(\x)$ and $\Delta(\x)$. These are defined pointwise in $\x$, as the smallest and largest (in absolute value) element of $\nabla\fo$, i.e.,
\(
\delta(\x)=1 \; \and \; \Delta(\x)=2\;\forall\x.
\)
Condition \qcc{Q2d.inf} now requires, for all $\x$ in $\D$, that
\[
\normi{\nabla\f(\x)}<\frac{\delta(\x)}{\Delta(\x)}=\frac{1}{2}.
\]
From the previous case we know that $\normi{\nabla\f(\x)} \le \max\{|a|,|b|\}$ for all $\x\in\D$, so the problem is guaranteed Fast-Lipschitz by qualifying condition \qc{Q2d}, provided that $-1/2<a,b\le0$.
Note that this would not have met the old case (ii) in \eqref{old ii.a}, since
\New{$\fox\neq c\1^T\x$.}

\paragraph*{$ab<0$}

When $a$ and $b$ have different signs, neither \qcc{Q1.pos}, nor \qcc{Q2.pos} holds. Instead, one can try qualifying condition \qc{Qinfd}, which does not place any sign restrictions on $\nabla\fx$. It is only required, in \qcc{Qinf.inf}, that
\[
\normi{\nabla\f(\x)} < \frac{\delta(\x)}{\delta(\x)+\Delta(\x)}.
 \]
These quantities are unchanged from the previous cases, so $\delta(\x)=1,\Delta(\x)=2$ and $\normi{\nabla\f(\x)} \le \max\{|a|,|b|\}$, so the problem is Fast-Lipschitz by qualifying condition \qc{Qinfd} if both $|a|$ and $|b|$ are less than 1/3.
Also in this case, the Old Qualifying Conditions would not have worked since  case (iii) in \eqref{old iii.a} requires a scalar objective function.

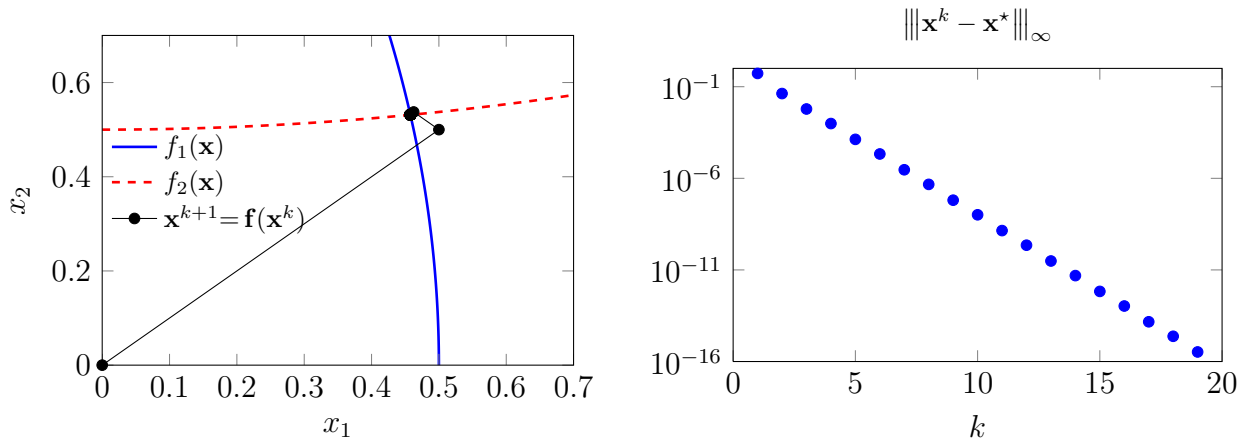
\begin{figure}[t]
        \centering
        \begin{subfigure}[b]{0.5\textwidth-1em}
            \setlength{\figwidth}{0.8\textwidth}
            \setlength{\figheight}{0.7\figwidth}
            \hspace{-10pt}
%
%
\begin{tikzpicture}

\begin{axis}[%
scale only axis,
width=\figwidth,
height=\figheight,
xmin=0, xmax=0.7,
ymin=0, ymax=0.7,
xlabel={$x_1$},
ylabel={$x_2$},
axis on top,
legend entries={$f_1(\x)$,$f_2(\x)$,$\x^{k+1}\hspace{-3pt}=\hspace{-1pt}\f(\x^k)$},
legend style={font=\small,at={(0.0,0.37)},fill=none,draw=none,anchor=south west,nodes=right}]
\addplot [
color=blue,
line width=1.0pt,
solid
]
coordinates{
 (0.5,0)(0.499985,0.010101)(0.499939,0.020202)(0.499862,0.030303)(0.499755,0.040404)(0.499617,0.0505051)(0.499449,0.0606061)(0.49925,0.0707071)(0.499021,0.0808081)(0.49876,0.0909091)(0.49847,0.10101)(0.498148,0.111111)(0.497796,0.121212)(0.497414,0.131313)(0.497,0.141414)(0.496556,0.151515)(0.496082,0.161616)(0.495577,0.171717)(0.495041,0.181818)(0.494475,0.191919)(0.493878,0.20202)(0.493251,0.212121)(0.492593,0.222222)(0.491904,0.232323)(0.491185,0.242424)(0.490435,0.252525)(0.489654,0.262626)(0.488843,0.272727)(0.488001,0.282828)(0.487129,0.292929)(0.486226,0.30303)(0.485292,0.313131)(0.484328,0.323232)(0.483333,0.333333)(0.482308,0.343434)(0.481252,0.353535)(0.480165,0.363636)(0.479048,0.373737)(0.4779,0.383838)(0.476722,0.393939)(0.475513,0.40404)(0.474273,0.414141)(0.473003,0.424242)(0.471702,0.434343)(0.47037,0.444444)(0.469008,0.454545)(0.467616,0.464646)(0.466192,0.474747)(0.464738,0.484848)(0.463254,0.494949)(0.461739,0.505051)(0.460193,0.515152)(0.458616,0.525253)(0.457009,0.535354)(0.455372,0.545455)(0.453704,0.555556)(0.452005,0.565657)(0.450275,0.575758)(0.448515,0.585859)(0.446725,0.59596)(0.444904,0.606061)(0.443052,0.616162)(0.441169,0.626263)(0.439256,0.636364)(0.437313,0.646465)(0.435338,0.656566)(0.433333,0.666667)(0.431298,0.676768)(0.429232,0.686869)(0.427135,0.69697)(0.425008,0.707071) 
};

\addplot [
color=red,
line width=1.0pt,
dashed
]
coordinates{
 (0,0.5)(0.010101,0.500015)(0.020202,0.500061)(0.030303,0.500138)(0.040404,0.500245)(0.0505051,0.500383)(0.0606061,0.500551)(0.0707071,0.50075)(0.0808081,0.500979)(0.0909091,0.50124)(0.10101,0.50153)(0.111111,0.501852)(0.121212,0.502204)(0.131313,0.502586)(0.141414,0.503)(0.151515,0.503444)(0.161616,0.503918)(0.171717,0.504423)(0.181818,0.504959)(0.191919,0.505525)(0.20202,0.506122)(0.212121,0.506749)(0.222222,0.507407)(0.232323,0.508096)(0.242424,0.508815)(0.252525,0.509565)(0.262626,0.510346)(0.272727,0.511157)(0.282828,0.511999)(0.292929,0.512871)(0.30303,0.513774)(0.313131,0.514708)(0.323232,0.515672)(0.333333,0.516667)(0.343434,0.517692)(0.353535,0.518748)(0.363636,0.519835)(0.373737,0.520952)(0.383838,0.5221)(0.393939,0.523278)(0.40404,0.524487)(0.414141,0.525727)(0.424242,0.526997)(0.434343,0.528298)(0.444444,0.52963)(0.454545,0.530992)(0.464646,0.532384)(0.474747,0.533808)(0.484848,0.535262)(0.494949,0.536746)(0.505051,0.538261)(0.515152,0.539807)(0.525253,0.541384)(0.535354,0.542991)(0.545455,0.544628)(0.555556,0.546296)(0.565657,0.547995)(0.575758,0.549725)(0.585859,0.551485)(0.59596,0.553275)(0.606061,0.555096)(0.616162,0.556948)(0.626263,0.558831)(0.636364,0.560744)(0.646465,0.562687)(0.656566,0.564662)(0.666667,0.566667)(0.676768,0.568702)(0.686869,0.570768)(0.69697,0.572865)(0.707071,0.574992) 
};

\addplot [
color=black,
solid,
mark=*,
mark options={solid,fill=black}
]
coordinates{
 (0,0)(0.5,0.5)(0.4625,0.5375)(0.456664,0.532086)(0.457533,0.531281)(0.457661,0.5314)(0.457642,0.531418)(0.457639,0.531415)(0.45764,0.531415)(0.45764,0.531415)(0.45764,0.531415)(0.45764,0.531415)(0.45764,0.531415)(0.45764,0.531415)(0.45764,0.531415)(0.45764,0.531415)(0.45764,0.531415)(0.45764,0.531415)(0.45764,0.531415)(0.45764,0.531415)(0.45764,0.531415)(0.45764,0.531415)(0.45764,0.531415)(0.45764,0.531415)(0.45764,0.531415)(0.45764,0.531415)(0.45764,0.531415)(0.45764,0.531415)(0.45764,0.531415)(0.45764,0.531415)(0.45764,0.531415)(0.45764,0.531415)(0.45764,0.531415)(0.45764,0.531415)(0.45764,0.531415)(0.45764,0.531415)(0.45764,0.531415)(0.45764,0.531415)(0.45764,0.531415)(0.45764,0.531415)(0.45764,0.531415)(0.45764,0.531415)(0.45764,0.531415)(0.45764,0.531415)(0.45764,0.531415)(0.45764,0.531415)(0.45764,0.531415)(0.45764,0.531415)(0.45764,0.531415)(0.45764,0.531415)(0.45764,0.531415)(0.45764,0.531415)(0.45764,0.531415)(0.45764,0.531415)(0.45764,0.531415)(0.45764,0.531415)(0.45764,0.531415)(0.45764,0.531415)(0.45764,0.531415)(0.45764,0.531415)(0.45764,0.531415)(0.45764,0.531415)(0.45764,0.531415)(0.45764,0.531415)(0.45764,0.531415)(0.45764,0.531415)(0.45764,0.531415)(0.45764,0.531415)(0.45764,0.531415)(0.45764,0.531415)(0.45764,0.531415)(0.45764,0.531415)(0.45764,0.531415)(0.45764,0.531415)(0.45764,0.531415)(0.45764,0.531415)(0.45764,0.531415)(0.45764,0.531415)(0.45764,0.531415)(0.45764,0.531415)(0.45764,0.531415)(0.45764,0.531415)(0.45764,0.531415)(0.45764,0.531415)(0.45764,0.531415)(0.45764,0.531415)(0.45764,0.531415)(0.45764,0.531415)(0.45764,0.531415)(0.45764,0.531415)(0.45764,0.531415)(0.45764,0.531415)(0.45764,0.531415)(0.45764,0.531415)(0.45764,0.531415)(0.45764,0.531415)(0.45764,0.531415)(0.45764,0.531415)(0.45764,0.531415)(0.45764,0.531415)(0.45764,0.531415) 
};

\end{axis}
\end{tikzpicture}
            \vspace{-10mm}
            \caption{Illustration of the feasible region of problem \eqref{eq toy no add} with $a=-0.3$ and $b=0.3$ (making the feasible region non-convex). The iterates \eqref{eq sol iter} of the solution quickly converges to $\xa=\f(\xa)$, where all constraints are active.}
            \label{fig1}
            \vspace{-5mm}
        \end{subfigure}%
        \hspace{1em}
        \begin{subfigure}[b]{0.5\textwidth-1em}
            \setlength{\figwidth}{65mm}
            \setlength{\figheight}{0.6\figwidth}
            \hspace{-2mm}
%
%
\begin{tikzpicture}

\begin{semilogyaxis}[%
scale only axis,
width=\figwidth,
height=\figheight,
xmin=0, xmax=20,
ymin=1e-016, ymax=1,
yminorticks=true,
xlabel={$k$},
title={\small$\normi{\x^k-\xa}$},
axis on top]
\addplot [
color=blue,
only marks,
mark=*,
mark options={solid}
]
coordinates{
 (1,0.531415)(2,0.0423603)(3,0.00608489)(4,0.000975634)(5,0.000133804)(6,2.13289e-005)(7,2.92836e-006)(8,4.66853e-007)(9,6.40952e-008)(10,1.02183e-008)(11,1.4029e-009)(12,2.23656e-010)(13,3.07062e-011)(14,4.89536e-012)(15,6.72018e-013)(16,1.07137e-013)(17,1.4766e-014)(18,2.38698e-015)(19,3.33067e-016)(20,0) 
};

\end{semilogyaxis}
\end{tikzpicture}
            \vspace{-10mm}
            \caption{Convergence of the iterates \eqref{eq sol iter}, measured in the $\left. \infty\text{-norm}\right.$. The convergence is geometric (linear in the log domain) and the optimal solution is found within an accuracy of $10^{-6}$  after 8 iterations.\\~}
            \label{fig2}
            \vspace{-5mm}
        \end{subfigure}
        \vspace{-4mm}
       \caption{Plots from the simple example in Section~\ref{sub ex simple}}
       \vspace{-6mm}
\end{figure}
\paragraph*{Solution of the problem}
If problem \eqref{eq toy no add} is Fast-Lipschitz by any of the cases above, the optimal point $\xa$ is found by solving $\xa=\f(\xa)$. We now solve the problem when $a=-0.3$ and $b=0.3$, by iterating
\begin{equation}
\label{eq sol iter}
\x^{k+1}:=\f(\x^k).
\end{equation}
This sequence will converge to $\xa=\f(\xa)$ since the qualifying conditions imply that $\f$ is contractive (Lemma~\ref{lem: f contrac}). The iterates $\x^k$ of \eqref{eq sol iter}, together with the feasible region of the problem is shown in Fig. \ref{fig1}. Clearly, $\xa=\fxa\in\cX$, so Proposition \ref{prop: additional constraint} applies and $\xa$ is optimal also for the original problem \eqref{eq toy example}. The convergence the iterations \eqref{eq sol iter} is shown in Fig. \ref{fig2}.

\New{
\subsection{Optimal control example}
\label{sub cex}
\newcommand\w{\mathbf{w}}
\newcommand\J{\mathbf{J}}
\renewcommand\F{\mathbf{F}}
\renewcommand\O{\mathbf{O}}
\newcommand\s{\mathbf{s}}
\newcommand\vs{\vec\s}
\newcommand\vu{\vec\u}
\newcommand\vw{\bar\w}
\newcommand\opt{{}^\star}
\newcommand\fy{\f_\y}
\newcommand\fyi{\f_ {\y\raisebox{.8ex}[0pt][0pt]{\ensuremath{_{[i]}}}}}

Consider a dynamical system with a state $\s\in \Rn$ and control variable $\u\in\Re^p$. The state evolves in discrete time and the state at time instance  $i+1$ is given by
\begin{equation}
\label{eq cex temp 7}
\s^{i+1}= \f(\s^i,\u^i)+\w^i,
\end{equation}
where $\f\,:\,\Rn\times\Re^p\to\Rn$ and $\w^i\in\Rn$ is an additive, bounded disturbance. The control variables are required to be positive and bounded, i.e., $\0\le\u\le\u_{\max}$. As a consequence of the bounded disturbance and control variables, the  state $\s^N$ after a finite number of $N$ iterations also remains bounded.
At each time instance, the system has a cost
$g(\s^i,\u^i)$ that is strictly increasing in all variables, i.e., all states represent something expensive and all controls are naturally associated with a positive cost.
 The design objective is to choose the control inputs $\{\u^i\}_{i=1}^{N}$ that minimizes the accumulated cost over $N$ periods.
 The optimal control problem, given the disturbances $\{\w^i\}_{i=1}^{N}$ and the initial state $\s_\text{init}$, becomes
\begin{equation}\label{eq cex temp 5}
\begin{array}{ll}
\min_{\{u^i\}_{i=1}^{N}}  &    \sum_{i=1}^{N} g(\s_i,\u_i)\\
\text{s.t.}    &  \s^1 = \s_\text{init} \\
                   & \s^{i+1}=\f(\s^i,\u^i)+\w^i, \quad i=1,\dots,N-1, \\
                   & \0\le \u^i \le \u_{\max}, \quad i=1,\dots,{N}.  \\
\end{array}
\end{equation}
Note that this can be seen as a centralized or as a distributed optimization problem.
Suppose that $\f(\s,\u)$ is increasing in $\u$. The optimal control would then trivially be $\u^i\opt=\0$ for all $i$, since  the cost $g(\s,\u^i)$ at time $i$ increases with $\u^i$ and the cost at time ${i+1}$ increases with $\s^{i+1}$, which in turn increases with $\u^i$. On the other hand, if $\f(\s,\u)$ is instead decreasing in $\u$, there would be a tradeoff between choosing a small $\u^i$ in order to make $g(\s^i,\u^i)$ small, or a large $\u^i$ to make $g(\s^{i+1},\u^{i+1})$ small. Throughout the rest of this example we will assume the non-trivial case when $\f(\s,\u)$ is decreasing in $\u$. With Fast-Lipschitz optimization it is possible to determine conditions on the costs and dynamics of problem~\eqref{eq cex temp 5}, under which the optimal control is simply given by $\u^i\opt = \0$ for all $i$.
This is a very useful result, which allows to avoid computing the control decision.
The following result is based on Proposition~\ref{prop: fewer constraints} and applies to dynamics given by a general function $\f(\s,\u)$:
\begin{result}
\label{result cex}
Consider problem~\eqref{eq cex temp 5} and assume that $g(\s,\u)$ is increasing in all variables. Assume further that the pair $\nabla_\s \f(\s,\u)$ and $\nabla_\s g(\s,\u)$ fulfill the \qc{Qk}  (in place of $\nabla\fx$ and $\nabla\fox$ respectively) and
\begin{equation}\label{eq cex result condition}
\frac{\max_{\s,\u} \normi{\nabla_\u\f(\s,\u)}}{1 - \max_{\s,\u} \normi{\nabla_\s\f(\s,\u)(\s,\u)}} < \frac{\min_{\s,\u} \min_i [\nabla_\u g(\s,\u)]_i}{\max_{\s,\u} \max_i [\nabla_\s g(\s,\u)]_i}
\end{equation}
for all allowed $\u$ and all reachable $\s$. Then, the optimal solution $\{\u^i\opt\}_{i=1}^N$ is given by $\u^i\opt = \0$ for all $i$, regardless of the problem horizon $N$, the initial state $\s_\text{init}$, and the disturbances $\{\w^i\}_{i=1}^N$.
\end{result}
The result can be derived be considering a Fast-Lipschitz problem equivalent to~\eqref{eq cex temp 5} as follows. Introduce the vectors
\[
\y = \begin{bmatrix} -\s^1 \\ \vdots \\ -\s^N \end{bmatrix} \in\Re^{nN},\quad
\z = \begin{bmatrix} -\u^1 \\ \vdots \\ -\u^N \end{bmatrix} \in\Re^{pN},\text{ and}\quad
\vw = \begin{bmatrix} \w^1 \\ \vdots \\ \w^N \end{bmatrix} \in\Re^{nN},
\]
and let $\y^i=-\s^i$, $\z^i=-\u^i$. Furthermore, let
$\x = \begin{bmatrix}\y^T & \z^T \end{bmatrix}^T$.
Problem~\eqref{eq cex temp 5} can then be transformed to the equivalent maximization problem
\begin{equation}\label{eq cex temp 6}
\begin{array}{ll}
\max_{\x}     &  f_0(\x) \\
\text{s.t.}    &  \y = \fy(\x),  \\
                   &\z_{\min} \le \z \le \0, \\
\end{array}
\end{equation}
where the $i$th $(n\times1)$-block of $\fy(\x)$ is given by
\[
\fyi (\x)= \begin{cases}
  -\s_\text{init} & i=1 \\
  -\f(-\y^i,-\z^i)-\w^i & i=2,\dots,N,
\end{cases}
\]
and $f_0(\x)=-\sum_{i=1}^{N} g(-\y_i,-\z_i)$.
By Corollary~\ref{prop: additional constraint} and Lemma~\ref{lem: All inequalities} we know that problem~\eqref{eq cex temp 6} is Fast-Lipschitz if the relaxed problem (obtained by replacing the equality constraints by inequality constraints, and holding the constraint $\z \ge \z_{\min}$ implicit)
\begin{equation}\label{eq cex temp 9}
\begin{array}{ll}
\max_{\x}     &  f_0(\x) \\
\text{s.t.}    &  \y \le \fy(\x),  \\
                   & \z \le \0, \\
\end{array}
\end{equation}
is Fast-Lipschitz. This is precisely a problem in the form of~\eqref{eq fewer constraints extended-1}, wherefore Proposition~\ref{prop: fewer constraints} applies.

To see this, we first denote the gradients of the original system dynamics~\eqref{eq cex temp 7} by
$\nabla_{\s} \f(\s,\u) \triangleq\A(\s,\u)$ and
$\nabla_{\u} \f(\s,\u) \triangleq \B(\s,\u).$
Simple calculations show that the gradient $\nabla_\y\fy(\x)$ consists of $N\times N$ blocks, each of dimension $n\times n$. All blocks are zero except the block sub-diagonal, which is given by the blocks
\[
\Big(\A(-\y^1,-\z^1) ,\dots,\A(-\y^{N-1},-\z^{N-1}) \Big),
\]
i.e.,
\begin{equation}\label{eq cex temp 8}
\nabla_{\y} \fy(\x) =
  \begin{bmatrix}
    \0  &  &  &  \\
    \A(-\y^1,-\z^1) & \ddots &  &  \\
     & \ddots & \ddots &  \\
     &  & \A(-\y^{N-1},-\z^{N-1})& \0  \\
  \end{bmatrix}.
\end{equation}
Similarly, the gradient $\nabla_\z\fy(\x)$ consists of $N\times N$ blocks of dimension $n\times p$, with
\[
\Big(\B(-\y^1,-\z^1) ,\dots,\B(-\y^{N-1},-\z^{N-1}) \Big)
\]
on the block sub-diagonal and zeros everywhere else.

In order to use Proposition~\ref{prop: fewer constraints} one must first verify that the subproblem $(f_{0 | \z}, \, \f_{\y | \z})$, obtained by fixing $\z$ in problem~\eqref{eq cex temp 9}, is Fast-Lipschitz for all permissible $\z$. This can be done by showing that $\nabla_{\y} \f_\y(\y,\z)$ and $\nabla_\y f_0(\y,\z)$ fulfills one of the qualifying conditions of Section~\ref{sec:Main}.
Next, one must verify part (b) of Proposition~\ref{prop: fewer constraints}. As we have restricted ourselves to the non-trivial case when $\f(\s,\u)$ is decreasing in $\u$, we have $\B(\s,\u)\le\0$. Consequently, $\nabla_\z\fy(\x)\not\ge\0$, and neither condition (b.i), nor condition (b.ii) applies. We therefore use condition (b.iii). Due to the block structure in~\eqref{eq cex temp 8}, we have
\[
\normi{\nabla_\y\fy(\x)} \le \max_{\s,\u} \normi{\A(\s,\u)} = \max_{\s,\u} \normi{\nabla_\s \f(\s,\u)}
\]
and
\[
\normi{\nabla_\z\fy(\x)} \le \max_{\s,\u} \normi{\B(\s,\u)} = \max_{\s,\u} \normi{\nabla_\u \f(\s,\u)}.
\]
The maximizations above are carried out over all permissible controls $\u$, and the corresponding achievable states $\s$.
Furthermore, we need
\(
  \Delta_{\y}(\x)=\max_{ij}\left[\nabla_{\y}f_0(\x)\right]_{ij}
\)
and
\(
\delta_{\z}(\x)=\min_{ij}\left[\nabla_{\z}f_0(\x)\right]_{ij}.
\)
Note that the gradients $\nabla_{\y}f_0(\x)$ and $\nabla_{\z}f_0(\x)$ are column vectors since $f_0$ is scalar. This gives
\[
\Delta_{\y}(\x)=\max_{i}\left[\nabla_{\y}f_0(\x)\right]_{i} \le \max_{\s,\u} \max_i [\nabla_\s g(\s,\u)]_i
\]
and
\[
\delta_{\z}(\x)=\min_{i}\left[\nabla_{\z}f_0(\x)\right]_{i} \le \min_{\s,\u} \min_i [\nabla_\u g(\s,\u)]_i,
\]
where the optimizations are over the permissible controls $\u$ and the reachable states $\s$.
Condition~(b.iii) is now fulfilled if
\[
\frac{\normi{\nabla_\z\fy(\x)}}{1 - \normi{\nabla_\y\fy(\x)}} < \frac{\delta_{\z}(\x)}{\Delta_{\y}(\x)},
\]
i.e., if
\begin{equation}\label{eq cex temp 10}
\frac{\max_{\s,\u} \normi{\nabla_\u \f(\s,\u)}}{1 - \max_{\s,\u} \normi{\nabla_\s \f(\s,\u)}} < \frac{\min_{\s,\u} \min_i [\nabla_\u g(\s,\u)]_i}{\max_{\s,\u} \max_i [\nabla_\s g(\s,\u)]_i}.
\end{equation}
When the inequality above holds, problem~\eqref{eq cex temp 9} fulfills Proposition~\ref{prop: fewer constraints}, whereby problems~\eqref{eq cex temp 9} and~\eqref{eq cex temp 6} are Fast-Lipschitz. This implies that the optimal solution is given by
\[
\begin{bmatrix}
\y\opt \\
\z\opt
\end{bmatrix}
=
\begin{bmatrix}
\f_\y(\y\opt,\z\opt) \\
\0
\end{bmatrix},
\]
i.e., the optimal solution of the original problem~\eqref{eq cex temp 5} is given by $\u^i\opt = -\z^i\opt = \0$ for all $i$. We will now apply Result~\ref{result cex} on two concrete examples of the dynamics $\nabla\f(\s,\u)$, one linear and one non-linear first order system.

\paragraph{Linear first order system}
Consider for illustrative purposes a linear first order system with a linear cost function. The optimal control problem~\eqref{eq cex temp 5} becomes
\begin{equation}
\label{eq cex lin prob}
\begin{array}{ll}
\min_{\{u^i\}_{i=1}^{N}}  &    \sum_{i=1}^{N} g(s_i,u_i)\\
\text{s.t.}    &  s^1 = s_\text{init} \\
                   & s^{i+1}=f(s^i,u^i)+w^i, \quad i=1,\dots,N-1, \\
                   & 0\le u^i \le u_{\max}, \quad i=1,\dots,{N},  \\
\end{array}
\end{equation}
where
\[
f( s ,u ) =  a s - b u
\quad
\text{and}
\quad
g(s,u) = c_s s + c_u u,
\]
with $a,b,c_s,c_u>0$. In this case, the gradients are given by $\nabla_s f(s,u) = a$, $\nabla_u f(s,u) = -b$, $\nabla_s g(s,u) = c_s$ and $\nabla_u g(s,u) = c_u$. Since the gradients are scalar and constant, it is straight forward to see that Result~\ref{result cex} applies to problem~\eqref{eq cex lin prob} provided that $a<1$ and
\[
\frac{b}{1-a} < \frac{c_u}{c_s}.
\]
If this is the case, the optimal solution must be $u^i\opt=0$ for all $i$. We remark that this is an important and non obvious result.

\paragraph{Non-linear first order system}
Consider again problem~\eqref{eq cex lin prob}, but with the dynamics given by
\begin{equation}\label{eq cex nonlin dyn}
f( s ,u ) =  \frac{a s}{1+s}s - b u.
\end{equation}
We continue to assume $a,b>0$, and the cost function $g(s,u)$ is left unchanged from the linear case.
The factor $as/(s+1)$ can be seen as a variable decay-rate, under which large states decay with a factor close to $a$ while small states decay almost instantly. Assume that the disturbances $w^i$ are non-negative, and let the upper limit of $u^i$ depend on the current state $s^i$,
\[
u \le u_{\max}(s^i) \quad \text{where} \quad  u_{\max}(s)= \frac{as^2}{b(1+s)}.
\]
Note that Corollary~\ref{prop: additional constraint} still applies when moving from problem~\eqref{eq cex temp 6} to problem~\eqref{eq cex temp 9} because $u^i\opt=0 \le u_{\max}(s^i)$, so the optimal values $u^i\opt=0$ of the relaxed problem~\eqref{eq cex temp 9} are always feasible in  problem~\eqref{eq cex temp 6}.
The modified constraints $u_{\max}(s^i)$ ensure  $f(s^i,u^i)\ge0$, whereby a non-negative $s_\text{init}\ge0$ implies that all future states $s^i$ are non-negative. All gradients from the linear example are left unchanged, except
\[
\nabla_s f(s,u) = a\frac{s^2+2s}{s^2+2s+1}.
\]
This gradient can be bounded by $0 \le \nabla_s f(s,u) \le a$, since the modified constraint $u_{\max}(s^i)$ ensures that only non-negative states $s$ can be reached.
Therefore, it holds that
\[
\max_{s,u} \normi{\nabla_s f(s,u)} \le a,
\]
whereby the conditions
\begin{equation}\label{eq cex nonlin cond}
|a|<1 \quad \text{and} \quad \frac{b}{1-a} < \frac{c_u}{c_s}
\end{equation}
remain unchanged from the linear example. When these conditions hold, we know that $u^i=0$ for all $i$ is the optimal solution. Because of the non-affine equality constraints in \eqref{eq cex nonlin dyn}, problem~\eqref{eq cex nonlin cond} is a non-convex problem in $N$ variables subject to $2N$ non-trivial constraints when $u_{\max}$ is a function of $s^i$. However, by using the Fast-Lipschitz properties of the problem we have solved it without performing any calculations, except for those involved in verifying the assumptions.

We now conclude this example by numerically solving two instances of problem~\eqref{eq cex lin prob}.

\paragraph{Numerical example}

\begin{table}[t]
\small
  \centering
  \New{
  \caption{\New{Shared problem parameter values}}\label{tab cex params}
  \begin{tabular}{c c c c c c}
    $N$ & $s_\text{init}$ & $\{w^i\}_{i=1}^N$ & $a$ & $c_s$ & $c_u$ \\
    \hline
    $20$ & $1$ & $w^i\sim\text{uni}[0,1]$ & $0.5$ & $3$ & $2$ \\
  \end{tabular}}
\end{table}

\begin{figure}[t]
  \centering
%
%
\begin{tikzpicture}

\begin{axis}[%
name=plot2,
scale only axis,
width=4.99714in,
height=1.03754in,
xmin=0, xmax=20,
ymin=0, ymax=1.5,
title={Resulting state},
xlabel={$i$},
ylabel={$s^i$},
axis on top]
\addplot [
color=blue,
solid,
mark=*,
mark options={solid}
]
coordinates{
 (1,1)(2,1.04473)(3,0.978263)(4,0.77824)(5,0.757091)(6,0.374054)(7,0.464829)(8,0.589572)(9,0.947671)(10,0.890463)(11,0.690893)(12,0.993913)(13,0.729248)(14,0.426063)(15,0.794843)(16,1.14547)(17,0.681405)(18,0.376673)(19,0.294236)(20,0.0523856) 
};

\addplot [
color=red,
dashed,
mark=asterisk,
mark options={solid}
]
coordinates{
 (1,1)(2,0.544732)(3,0.535047)(4,0.4143)(5,0.456191)(6,0.282404)(7,0.413964)(8,0.463062)(9,0.680084)(10,0.457515)(11,0.467028)(12,0.69359)(13,0.443544)(14,0.340438)(15,0.604208)(16,0.78115)(17,0.463802)(18,0.312077)(19,0.279819)(20,0.049529) 
};

\end{axis}

\begin{axis}[%
at=(plot2.above north west), anchor=below south west,
scale only axis,
width=4.99714in,
height=1.03754in,
xmin=0, xmax=20,
ymin=-5e-2, ymax=1,
title={Optimal control},
xlabel={$i$},
ylabel={$u^i$},
axis on top]
\addplot [
color=blue,
solid,
mark=*,
mark options={solid}
]
coordinates{
 (1,0)(2,0)(3,0)(4,0)(5,0)(6,0)(7,0)(8,0)(9,0)(10,0)(11,0)(12,0)(13,0)(14,0)(15,0)(16,0)(17,0)(18,0)(19,0)(20,0) 
};

\addplot [
color=red,
dashed,
mark=asterisk,
mark options={solid}
]
coordinates{
 (1,1)(2,0.544732)(3,0.430615)(4,0.382571)(5,-4.3658e-21)(6,0.062091)(7,0.226713)(8,0.463062)(9,0.680083)(10,0.171911)(11,0.467028)(12,0.360019)(13,-1.15031e-20)(14,0.340437)(15,0.604208)(16,0.166226)(17,0)(18,0)(19,0)(20,0) 
};

\end{axis}

\begin{axis}[%
at=(plot2.below south west), anchor=above north west,
scale only axis,
width=4.99714in,
height=1.03754in,
xmin=0, xmax=20,
ymin=0, ymax=6,
title={Resulting cost per stage},
xlabel={$i$},
ylabel={$ g(s^i,u^i)$},
axis on top]
\addplot [
color=blue,
solid,
mark=*,
mark options={solid}
]
coordinates{
 (1,3)(2,3.13419)(3,2.93479)(4,2.33472)(5,2.27127)(6,1.12216)(7,1.39449)(8,1.76872)(9,2.84301)(10,2.67139)(11,2.07268)(12,2.98174)(13,2.18774)(14,1.27819)(15,2.38453)(16,3.4364)(17,2.04422)(18,1.13002)(19,0.882707)(20,0.157157) 
};

\addplot [
color=red,
dashed,
mark=asterisk,
mark options={solid}
]
coordinates{
 (1,5)(2,2.72366)(3,2.46637)(4,2.00804)(5,1.36857)(6,0.971395)(7,1.69532)(8,2.31531)(9,3.40042)(10,1.71637)(11,2.33514)(12,2.80081)(13,1.33063)(14,1.70219)(15,3.02104)(16,2.6759)(17,1.39141)(18,0.936232)(19,0.839456)(20,0.148587) 
};

\end{axis}
\end{tikzpicture}\\
  \caption{\New{The plots show the optimal control variables, and the resulting states and costs. The first system ($b=0.3$) is marked with a solid blue line and dots. As expected, the optimal control is zero at all times. The second system ($b=0.5$) is marked by a dashed red line and stars. The total accumulated cost (i.e., the sum over each stage-cost line) is $42.03$ and $40.85$ respectively. The second system $(b=0.5)$ has a lower total cost and this is expected, since the second system has a more powerful actuator (at the same actuation cost).}}
  \label{fig cex1}
\end{figure}
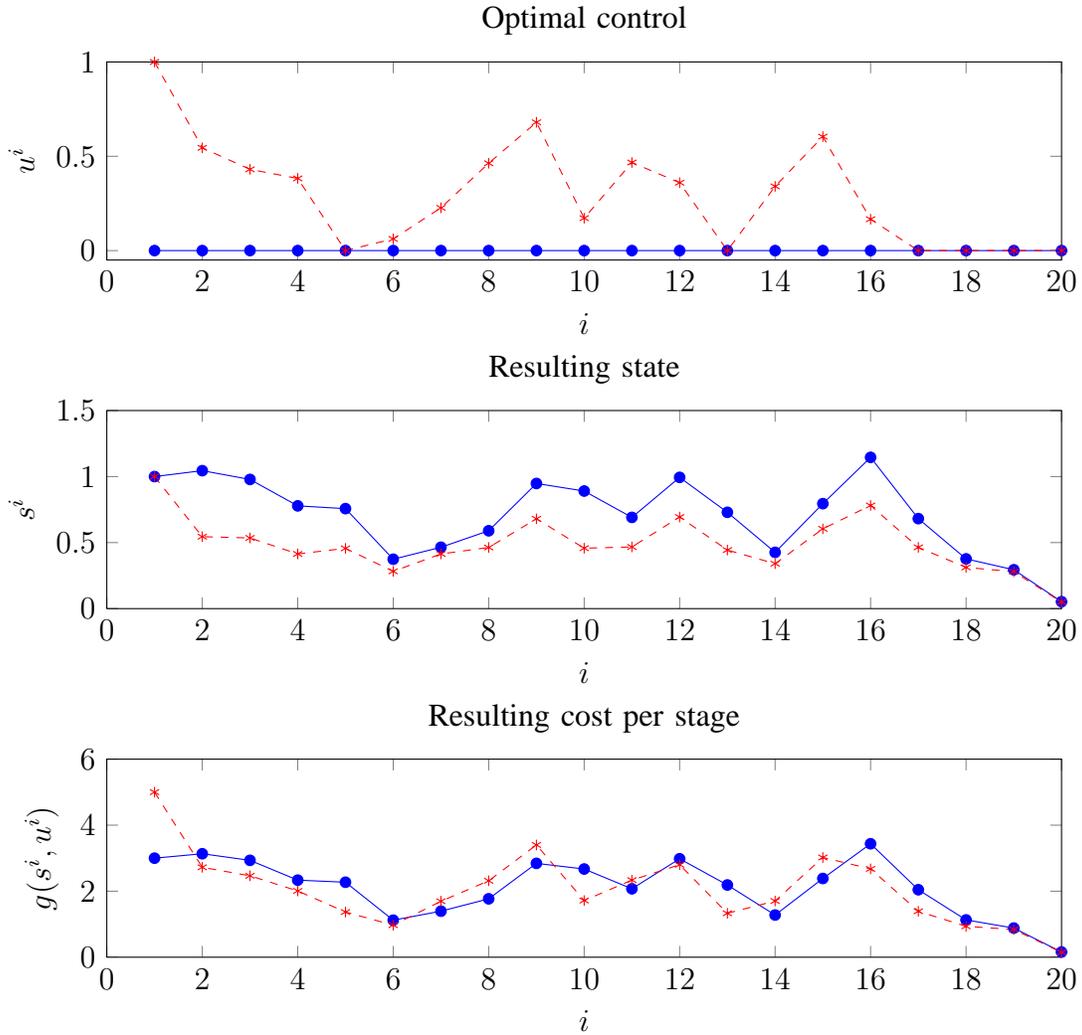
Consider two different instances of problem~\eqref{eq cex lin prob}, with the non-linear dynamics given by~\eqref{eq cex nonlin dyn}. The two problems share all parameter values, given in Tab.~\ref{tab cex params}, except for that of $b$; the first problem instance has $b=0.3$ and the second instance has $b=0.5$. Both instances fulfill $0\le a<1$, but only the first instance fulfill
\[
\frac{b}{1-a}=\frac{0.3}{1-0.5} = 0.6 <  \frac{2}{3} = \frac{c_u}{c_s}.
\]
We can therefore guarantee that $u^i=0$ for all $i$ solves the first instance, but we cannot say anything about the second instance. Remember, the conditions in this paper are only sufficient.
The two problems are solved numerically with Matlab's built-in solver \texttt{fmincon}, and the results are shown in Fig.~\ref{fig cex1}. As expected by Result~\ref{result cex}, the optimal actions for the first system ($b=0.3$) is to always keep the control variables at the lower bound. The second system ($b=0.5$) illustrates a case when Result~\ref{result cex} does not apply.
}

\section{Conclusions and future work}\label{sec:Conclusions-and-future}
\New{
In this paper we significantly extended the previous Fast-Lipschitz framework proposed in \cite{Fischione11}. A new set of qualifying conditions, that unify and generalize previous conditions, was presented.
Furthermore, we investigated problems that deviate from the required Fast-Lipschitz form, either by not having the same number of constraints as variables or by having variables that do not affect the objective function. For these cases we established conditions for which they are still Fast-Lipschitz.
Based on these new results, a larger set of convex and non-convex optimization problems can be solved by the Fast-Lipschitz method of this paper, both in a centralized and in a distributed set-up. This avoids using Lagrangian methods, which are inefficient in terms of computations and communication complexity especially when used over networks.

Several possible extensions remain to consider, for example:
\begin{itemize}
  \item We believe that the Fast-Lipschitz optimization framework can be extended to cover also non-smooth problems.  A potential benefit of such an extension would be to form expressions such as
      $\x\le\fx=\min\{\f_1\ox,\f_2\ox\},$
  whereby problems with more constraints than variables ($\fx:\Rn\to\Re^p$ where $p>n$) could be considered.

  \item So far, only problems with $\nabla\fox\ge\0$ have been considered. The standard inequality we have used is a partial ordering induced by the non-negative orthant $\Rmp$, i.e., $\nabla\fox\succeq_{\Rmp}\0$. A possible extension is to allow for problems where $\nabla\fox\succeq_{\mathcal{K}}\0$ for a more general cone~$\mathcal{K}$. This would allow a tradeoff between conditions for $\nabla\fox$ and conditions on $\nabla\fx$, which may give more flexibility in the qualifying conditions.

    We have already seen a similar example of this, when case (ii) of the old conditions is generalized to \qc{Q2}. In this case, condition \eqref{old ii.a} required $\fox=c\1^T\x$ for some $c>0$, which can also be stated
    \[
    \fox\in\Re,\quad \text{and}\quad\nabla\fox\succeq_{\mathcal{K}_\1}\0
    \]
  where  $\mathcal{K}_\1$ is the cone (ray) generated by the single vector $\1$. In \qc{Q2} it is instead required that $\nabla\fox>\0$, which is more general. This, however, comes at the price of a less general requirement on $\fx$, since
  \qcc{Q2.inf} requires
  \[
  \normi{\nabla\fx}<\q\ox\le 1,
  \]
  while \eqref{old ii.c} only requires
  $\normi{\nabla\fx}<1.$

\end{itemize}
}
\bibliographystyle{ieeetr}
\bibliography{references_ieee_formated}

\New{
\appendix
\begin{lemma}
Let $\norm{\cdot}_a$ be a matrix norm. Then $\norm{\cdot}_b$, defined as
  $\norm{\A}_b=\norm{\A^T}_a$
  is also a matrix norm.
\end{lemma}
\begin{TempProof}
A function $\norm\cdot$ is a matrix norm if it fulfills
\begin{enumerate}
  \item $\norm{\A}\ge0$
  \item $\norm{\A}=0 \iff \A=\0$
  \item $\norm{c\A}=c\norm{A}$ for all scalars $c$
  \item $\norm{\A+\B}\le\norm{\A}+\norm{\B}$
  \item $\norm{\A\B} \le \norm{\A}\norm{\B}$
\end{enumerate}
  \newcommand\markdef[1]{\overset{\star a}{#1}}
The norm $\norm{\A}_b = \norm{\A^T}_a$ fulfills all the matrix norm properties from above, since
  \begin{enumerate}
     \item $\norm{\A}_b=\norm{\A^T}_a \markdef\ge0$
  \item $\norm{\A}_b=\norm{\A^T}_a = 0 \markdef\iff \A^T=\0 \iff \A=\0$
  \item $\norm{c\A}_b=\norm{c\A^T}_a \markdef= c\norm{\A^T}_a = c\norm{\A}_b$ for all scalars $c$
  \item $\norm{\A+\B}_b=\norm{\A^T+\B^T}_a \markdef\le \norm{\A^T}_a+\norm{\B^T}_a=\norm{\A}_b+\norm{\B}_b$
  \item $\norm{\A\B}_b=\norm{\B^T\A^T}_a \markdef\le \norm{\B^T}_a\norm{\A^T}_a=\norm{\A}_b\norm{\B}_b$
  \end{enumerate}
  Each of the chains above uses only the definition of $\norm\cdot_b$, and the corresponding property that holds for the matrix norm $\norm\cdot_a$ (marked by  $\markdef~$).
\end{TempProof}
}
\end{document}